\documentclass[12pt]{amsart}
\usepackage{geometry}                
\usepackage{graphicx}
\usepackage{amssymb}
\usepackage{epstopdf}
\title{\bf Quasilinear P.D.E.s,  Interpolation spaces and  H\"olderian mappings }
\author[I. Ahmed,  A. Fiorenza, M.R. Formica,  A. Gogatishvili, A. El Hamidi, J. M. Rakotoson]{I. Ahmed{$^{(1)}$},\ \ \  A. Fiorenza$^{(2)}$,\ \ \ M.R. Formica$^{(3)}$,\\ A. Gogatishvili$^{(4)}$,\ \  \ A. El Hamidi$^{(5)}$,\ \ \ J. M. Rakotoson$^{(6)}$}
\address{$^{1}${\small  Department of Mathematics, Sukkur IBA University, Sukkur, }
{\bf Pakistan} }
\email{
{\small irshaad.ahmed@iba-suk.edu.pk}}
\address{{$^{2}$ Universit\`a di Napoli "Federico II"}\\
{\small  via Monteoliveto, 3, I-80134 Napoli,  {\bf Italy}}\\
{\small and Istituto per le Applicazioni del Calcolo "Mauro Picone" }\\
{\small Consiglio Nazionale delle Ricerche} \\
{\small via Pietro Castellino, 111 I-80131Napoli, {\bf Italy}}\\ }
\email{
{\small    fiorenza@unina.it}}
\address{$^{3}${\small Universit\`a degli Studi di Napoli "Parthenope", via  Generale Parisi 13, 80132, Napoli, {\bf Italy}}\\ }
\email{
{\small   mara.formica@uniparthenope.it}}
\address{$^{4}${\small  Institute of  Mathematics of the   Czech Academy of Sciences - }
 {\small \v Zitn\'a, 115 67 Prague 1, {\bf Czech Republic}}}
\email{
{\small  gogatish@math.cas.cz}}
\address{$^{5}${\small D\'epartement de Math\'ematiques et Laboratoire LaSIE, }
{\small Universit\'e de La Rochelle, Av. Michel Cr\'epeau  17042 La Rochelle cedex 1 {\bf France}}}
\email{
{\small aelhamid@univ-lr.fr}}
\address{$^{6} ${\it corresponding  author,  }
{\small Laboratoire de Math\'ematiques et Applications - Universit\'e de Poitiers,}
{\small 11  Bd  Marie et Pierre Curie,T\'el\'eport 2, }
{\small 86073 Poitiers Cedex 9, {\bf France}}
} 
\email{
{\small  rako@math.univ-poitiers.fr,   jean.michel.rakotoson@univ-poitiers.fr }}

\date{}                                           
\def\a{\alpha}
\def\b{\beta}
\def\eps{\varepsilon}
\def\f{\varphi}
\def\l{\lambda}
\def\p{\partial}
\def\s{\sigma}
\def\tpl{{\theta,p;\lambda}}
\def\tpala{{\theta,\frac p\alpha;\lambda\alpha}}
\def\O{\Omega}
\def\OV{\overline}
\def\R{{\rm I\!R}}
\def\S{\mathbb S}
\def\Log{{\rm Log\,}}
\def\Max{\mathop{\rm Max\,}}
\def\Min{\mathop{\rm Min\,}}
\def\WH{\widehat}
\def\WT{\widetilde}
\def\div{{\rm div\,}}

\def\LEQ{\leqslant}
\def\GEQ{\geqslant}
\def\DST{\displaystyle}
\def\esssup{\mathop{\rm ess\,sup\ }}
\def\Inf{\mathop{{\rm Inf\,}}}
\def\dist{{\rm dist\,}}
\def\Nu{\nabla u}
 \def\Nv{\nabla v}
 \def\Numv{\nabla(u-v)}
 \def\with{\hbox{ with }}
\def\pep{{(p^*)'}}
\def\Pr{{\it Proof: }}
\def\HF{\hfill{$\diamondsuit$}\\ }
\newtheorem{theo}{Theorem}[section]
\newtheorem{lem}{Lemma}[section]

\newtheorem{propo}{Proposition}[section]

\newtheorem{coro}{Corollary}[theo]
\newtheorem{corol}{Corollary}[lem]
\newtheorem{corop}{Corollary}[propo]
\newtheorem{defi}{Definition}[section]
\newtheorem{rem}{Remark}[section]
\def\Zue{{(L^1,L^m)_{\theta,p_2;\lambda}}}
\def\Zde{{(L^1,L^{n,1})_{\theta,p_2;\lambda}}}
\def\Zte{{(L^{s,\infty},L^m)_{\theta,p_2;\lambda}}}
\def\Zqe{{(L^{s,\infty},L^\infty)_{\theta,p_2;\lambda}}}
\def\tpdl{{\theta,p_2;\lambda}}

\def\vecp{{\vec p}}

\def\nppp{{\frac{n'}{p'}}}

\def\Lpep{{L^{(p^*)'}}}
\def\meas{{\rm meas\,}}
\def\Min{\mathop{{\rm Min\,}}}

\def\pp{{p(\cdot)}}
\def\sign{{\rm sign}}
\def\qq{{[p^*(\cdot)]'}}
\def\qp{q(\cdot)}
\def\C{\centerline}
\def\point{{p(\cdot)}}
\def\pointetoile{{p^*(\cdot)}}
\def\psimT{\hbox{$\eta B\big(T_{m+1}(u)-\varphi\big)$}}
\def\psimu{\hbox{$\eta B\big(u^{m+1}-\varphi\big)$}}
\def\psiuf{\hbox{$\eta B\big(u -\varphi\big)$}}
\def\si{{\hbox{ if }}}
\def\ewe{{\hbox{ elsewhere }}}
\def\aein{{\hbox{ a.e in }}}
\def\mpu{{m+1}}

\begin{document}
\maketitle
\begin{abstract} \ \\
As in the work of Tartar ([59]), we develop here some new results on nonlinear interpolation of $\alpha$-H\"olderian mappings between normed spaces, by studying the action of the mappings on $K$-functionals and between interpolation spaces with logarithm functions. We apply these results to obtain some regularity results on the gradient of the solutions to quasilinear equations of the form $$-\div(\widehat a(\nabla u ))+V(u)=f, $$ where $V$ is a nonlinear potential and  $f$ belongs to non-standard spaces like Lorentz-Zygmund spaces. We show several results; for instance, that the mapping ${\mathcal T}: \ {\mathcal T}f=\nabla u$ is locally or globally $\alpha$-H\"olderian under suitable values of $\alpha$ and appropriate hypotheses on $V$ and $\widehat a.$
\end{abstract}

\tableofcontents
\section{\bf\large Introduction - Notation - Preliminary results}
\subsection{\bf Introduction}\ \\
The Marcinkiewicz interpolation theorems for linear operators acting on Lebesgue spaces turned out to be an essential tool for studying  regularity of solutions for linear partial differential equations (P.D.E.s) in $L^p$-spaces.
Then, Jaak Peetre (\cite{Peetre1, Peetre1b}) introduced a method $(K$-method) to give a general definition of interpolation spaces between two normed spaces embedded in a same topological space. His definition allows to extend the Marcinkiewicz's results of linear operators to those ones acting on abstract normed spaces. But  his results allow also to go further in the study of  regularity of   solutions of linear equations on spaces different from  $L^p$ spaces. The main problem to apply Peetre's definition is the identification of the interpolated spaces. Some results in this direction exist: for instance, we did such a  study with applications to linear P.D.E.s in recent papers (see \cite{FFGKR}, \cite{AFFGR} or \cite{FFR}) using new spaces as grand or small Lebesgue spaces, sometimes combining the regularity method with a duality method.\\
Later, in our knowledge, L. Tartar \cite{Tartar3}, under the supervision of J.L. Lions, was the first to give some interpolation results on nonlinear H\"olderian mappings (which include Lipschitz mappings) and he applied them to a variety of boundary value problems as bilinear applications, to semi-linear P.D.E.s but also to variational inequalities .\\
This last paper of L. Tartar, recent results development concerning the interpolation spaces with logarithm functions  (see, for instance \cite{GOT}, and the previous references) and the appearance of the new operators in P.D.E.s as anisotropic $\vec p$-Laplacian or variable exponents $p(\cdot)$-Laplacian, were the main motivations which lead us to reconsider the work of L. Tartar \cite{Tartar3} and to show that we may have H\"older mappings associated to quasilinear equations in order to obtain new regularity results. \\
So, we extend first Tartar's results on nonlinear interpolations mappings ${\mathcal T}$ to couples of spaces with a logarithm function by studying the action of the mapping ${\mathcal T}$ on the  $K$-functional associated to those couples. This is the purpose of the second section. Here is an example of such result:

{\it Let $X_1\subset X_0, \ Y_1\subset Y_0$, be four normed  spaces, and let $0<\a\LEQ1$. Assume that ${\mathcal T}:X_i\to Y_i$ is globally $\a$-H\"olderian for $i=0,1$ with H\"older constant $M_i$, i.e. 
$$\exists \, M_i>0 \,\, {\it such} \,\, {\it that} \,\, ||{\mathcal T}a-{\mathcal T}b||_{Y_i}\LEQ M_i||a-b||^\a_{X_i}, \qquad i=0,1.$$
Then, for all $a\in X_0, \ b\in X_1$ one has
$$K({\mathcal T}a-{\mathcal T}b;t^\a)\LEQ 2\Max(M_0;M_1) K(a-b;t)^\a.$$
} \\
As a consequence, we derive the following result:\\
\ \\
{\it
Let $X_1\subset X_0, \ Y_1\subset Y_0$ four normed  spaces. Assume that ${\mathcal T}:X_i\to Y_i$ is globally
$\a$-H\"olderian for $i=0,1$.}
{\it For $0\LEQ\theta\LEQ1,\ 1\LEQ p\LEQ+\infty$, if $X_1$ is dense in $X_0$, then}
$$\hbox{\it ${\mathcal T}$ is an $\a$-H\"olderian mapping from $(X_0,X_1)_{\theta,p;\lambda}$ into $(Y_0,Y_1)_{\theta,\frac p\a;\lambda \a}$}.$$

The last part of the second section is devoted to some identification of interpolation spaces using  couples of  Lebesgue or Lorentz spaces. This  allows us to recover spaces as Lorentz-Zygmund spaces or $G\Gamma$-gamma spaces. The list is not exhaustive but was chosen to be applied later on, in the fourth and the fifth sections.\\
To define the appropriate mappings in those last sections,  we   consider two types of formulations, the usual weak formulation and the entropic-renormalized formulation for the quasilinear P.D.E.s of the form $Au+V(u)=f,\ f\ $ in$\ L^1(\O)$, where $A$ is a Leray-Lions type operator, $V$ a potential, and we may prove  the existence and uniqueness of solution according to the space where the data $f$ belongs. We can define a non-linear operator, ${\mathcal T}:L^1(\O)\to Y_{0i}$, $i=1,\cdots,n$~: to $f\in L^1(\O)$ we associate the $i$-th component of the {\bf gradient} of the solution in an appropriate space $Y_{0i}$.
The main step is to prove that such a nonlinear operator is a H\"olderian mapping. This is done in each application from section four to six. The fundamental lemma (see {\bf Lemma \ref{lBt} } below) to obtain such a result
in Marcinkiewicz spaces  for $L^1$ data reads as follows:  \\
{\it Let $\nu $ be a non negative Borel measure and $h:\O\to\R_+,\ g:\O\to\R_+$, be two $\nu$-measurable functions. Then\\
\C{$\DST\nu\Big\{h>\lambda\Big\}\LEQ\frac1\lambda \int_{\{g\LEQ k\}}hd\nu+\nu\{g>k\}\qquad \forall\,\lambda>0,\ \forall\,k>0.$}
}

Replacing $L^1(\O)$ by other  $L^r$-spaces  we can have more regularity on the gradient of the solution.\\
We then apply the abstract results on interpolation mappings obtained in the second section. Let us notice that our estimates   are optimal in many cases. Therefore we improve some well-known regularity results as in Lorentz spaces but also we have an easy tool to derive regularity of the gradient when the data $f$ is  in spaces as $L^{m,r}\big(\Log\,L\big)^\a,\ m\GEQ 1$ or in small spaces $L^{(r,\theta}(\O)$ or Orlicz spaces.\\
For convenience, we took only models for the nonlinear operator $A$.
 More precisely, we  study the regularity of the weak or entropic-renormalized solution of a p-laplacian type operators such as
$-\div\Big(|\nabla u|^{p-2}\nabla u\Big)+V(x;u)=f$,
or its anisotropic version  in a bounded smooth domain $\O$ of $\R^n$,
$$-\sum_{i=1}^n\dfrac\partial{\partial x_i}\left[\left|\dfrac{\partial u}{\partial x_i}\right|^{p_i-2}\dfrac{\partial u}{\partial x_i}\right]+V(x;u)=f, \qquad 1\!<p_i,\,p\!<\!+\infty\, ,\,\, i\!=\!1,\!\ldots\!,\!n\, ,$$
or the variable exponents version of $p(\cdot)$-Laplacian, where $V$ is nonlinear. We only consider the Dirichlet  homogeneous condition on the boundary $u=0$.
\\
An example of regularity that we can prove (it will be a consequence of  Proposition \ref{p3}, see below) is:
{\it If $u$ is a  solution of the quasilinear equation (\ref{eq15}) (see below), $2\LEQ p < n$ and  $f\in L^{k,r}(\O)$, then the gradient of the solution $u$ belongs to $[L^{k^*(p-1),r(p-1)}(\O)]^n$, with $k \LEQ (p^*)'$ (here $(p^*)'$ denotes the conjugate of the Sobolev exponent of p, and $k^*$ denotes the Sobolev exponent of $k$). Moreover, we have
$$||\nabla u||_{L^{k^*(p-1),r(p-1)}}\LEQ c||f||_{L^{k,r}}^{\frac1{p-1}}.$$
}
An example of non-standard regularity result that can be obtained from Theorem \ref{t6} (see below)
for the solution $u$  is:\\
{\it $$\left[ \int_0^1\left(\left(\int_t^1 |\nabla u|_\ast(s)^{p}ds\right)^{\frac{1}{p}}
			 (1-\log t)^{\lambda\alpha}\right)^{\frac{p}{\alpha}} \dfrac{dt}{t}\right]^{\frac{\alpha}{p}} \LEQ $$
$$c\left[ \int_0^1\left(\left(\int_t^1 f_\ast(s)^{(p^*)'}ds\right)^{\frac{1}{(p^*)'}}
			 (1-\log t)^{\lambda }\right)^p \dfrac{dt}{t}\right]^{\frac \alpha p},$$
			 whenever the right hand side of the inequality is finite. Here $\alpha = \frac 1{p-1},\; 2 \LEQ p < n,\; \lambda\in\R$ .\\ Moreover,
	if $ f \in L^{\frac{m'}{m'-\theta},p_2}\Big(\Log L\Big)^\lambda$,  then   $|\nabla u| \in L^{p_\theta,p_2(p-1)}\Big(\Log L\Big)^{\frac\lambda{p-1}}$.  }\\

From Section 4 to Section 6, we give some applications of the  abstract results obtained  in Sections 2 and 3.
For instance, here is the basis of the existence of an  H\"olderian mapping result for anisotropic equation:
{\it Let $u$ be the entropic-renormalized solution of equation (\ref{eq1150}) (see below).
 Then there exists a constant $c>0$ independent of $u$ and $f$ such that
\begin{enumerate}
\item ${\rm meas\,}\{|u|>k\}\LEQ c||f||_{L^1(\O)}^{\frac{p^*}p}k^{-\frac{p^*}{p'}}$,\ \ $\forall\,k>0$.\\
\item $\DST\left\|\dfrac{\p u}{\p x_i}\right\|_{L^{\frac{n'p_i,}{p'}\infty}(\O)}\LEQ c||f||_{L^1(\O)}^{\frac{p'}{p_i}},\qquad\qquad i=1,\ldots,n.$
\end{enumerate} }

For the sake of completeness, although the existence and uniqueness for  quasilinear equations are widely done in the literature and are  not the main issue of our work, we shall give some examples of proofs of  uniqueness and existence. Namely, when the operator $A$ has variable exponents, we have new  results   and we show in particular that:\\
{\it There exists a constant $c>0$ depending only on $p,\ n,\ \O$ such that
$$\meas\Big\{|\Nu|^\point>\lambda\Big\}\LEQ c\,\psi_1(||f||_1)^{\frac1{1+|a_1|}}\lambda^{-\frac{|a_1|}{1+|a_1|}}\quad\forall\,\lambda>0.$$}
Such topic is developed in the last section 6. The method is widely inspired by the previous works (see for instance \cite{Benilanetal}, \cite{Rako2}, \cite{Rako3}, \cite{Rako5}, \cite{Rako6}), and uses recent theorems as the one given in \cite{FGNR}.
Moreover, the same method can be used  to prove the existence and uniqueness of entropic-renormalized solution for general operators including the anisotropic case.\\
For other results concerning interpolation of Lipschitz operators and other applications of Interpolation theory, also in P.D.E.s, see \cite{Cianchi-Mazja2, Malig1, Malig2, Malig3}. 

 \subsection{\bf Notations -Preliminary results}\ \\
 We shall adopt our usual notations.
 For a measurable $f\in\O\to\R$, we set for $t\GEQ 0$
 $$D_f(t)={{\rm meas\,}}\Big\{ x\in\O:|f(x)|\GEQ t\Big\},$$
  and $f_*$, the decreasing rearrangement of $|f|$, is defined as follows: for $s\in(0,|\O|),$ $|\O|$ being the measure of $\O$,
  $$f_*(s)=\inf\Big\{t:D_f(t)\LEQ s\Big\}.$$
  
  We also set 
 $$ f_{**}(s)=\dfrac 1s\int_0^s f_*(t)dt. $$
 The Lorentz space $L^{p,q}(\O),$ \ $1<p<+\infty$, \ $1\LEQ q\LEQ+\infty,$ is defined as the set of measurable functions $f$ for which
  $$||f||_{p,q}=\begin{cases}\DST
  \left[\int_0^{|\O|}[t^{\frac1p}f_{**}(t)]^q\dfrac{dt}t\right]^{\frac1p}&\hbox{if $q<+\infty$},\\
  \DST\sup_{0<t<|\O|}t^{\frac1p}f_{**}(t)&\hbox{if $q=+\infty,$}\end{cases}\qquad\hbox{is finite,}$$
   while $||v||_q$  denotes the norm in $L^q(\O),\ 1\LEQ q \LEQ +\infty$. \\
  If $A_1$ and $A_2$ are two quantities depending on some parameters, we shall write
  $$A_1\lesssim A_2$$
$\hbox{if there exists $c>0$ independent of the parameters such that }A_1\LEQ c A_2$, and
  $$A_1\simeq A_2$$
$\hbox{if and only if $A_1\lesssim A_2$ and $A_2\lesssim A_1.$}$

  For the anisotropic problem, we will need the following Troisi's Sobolev inequalities \cite{V, Tr}.
  Setting $$\dfrac1p =\dfrac 1n\DST\sum_{i=1}^n\dfrac1{p_i}\quad\hbox{ and }\quad p^*=\dfrac{np}{n-p}\hbox{\ \  if \ \ }\DST\sum_{i=1}^n\dfrac1{p_i}>1,\quad \vec p=(p_1,\ldots,p_n)\, ,$$
  we have\\
\begin{theo}\label{tlN1}{\bf (Poincar\'e-Sobolev inequality for anisotropic Sobolev spaces)}\ \\
If $1\LEQ p<n$, $1\LEQ p_i<n$ $(i=1,\ldots,n)$, then the following inequalities hold true.
\begin{enumerate}
\item There exists a constant $C=C(n,\vec p)$ such that
\begin{equation}
\label{eqn1}
\left(\int_{\R^n}|u|^{p^*}dx\right)^{\frac n{p^*}}\LEQ C\prod_{i=1}^n\left(\int_{\R^n}
|\p_{x_i}u|^{p_i}dx\right)^{\frac1{p_i}}\qquad \forall\,u\in C^\infty_c(\R^n).
\end{equation}
\item For any $\vec\theta=(\theta_1,\ldots,\theta_n)$ such that $\theta_i>0$ for all $i=1,\ldots,n$ and $\DST\sum_{i=1}^n\dfrac1{\theta_i}=\dfrac np$, there exists a constant $C_{\vec\theta}=C(n,\vec p,\vec\theta)$ such that
\begin{equation}
\label{eqn2}
\left(\int_{\R^n}|u|^{p^*}dx\right)^{\frac p{p^*}}\LEQ C_{\vec \theta}\sum_{i=1}^n\left(\int_{\R^n}|\p_{x_i}u|^{p_i}dx\right)^{\frac{\theta_i}{p_i}} \quad\forall\,u\in C^\infty_c(\R^n).
\end{equation}
In particular, we shall use the case $\theta_i=p_i$ for all $i=1,\ldots,n$.
\end{enumerate}
\end{theo}

We shall denote by $W^{1,\vec p}_0(\O)$ the closure of $C_c^\infty(\O)$ with respect to the norm:
$$||v||_{1,\vec p}=\sum_{i=0}^n\left\|\dfrac{\p v}{\p x_i}\right\|_{p_i}. $$
The  following Poincar\'e-Sobolev  inequality holds true in
$W^{1,\vec p}_0(\O)$.
\begin{coro}[\bf of Theorem \ref{tlN1}]\label{c1tN1}
\hfill
\begin{enumerate}
\item There exists a constant $C=C(n,\vec p)$ such that
\begin{equation*}
\left[\int_\O|v|^{p^*}(x)dx\right]^{\frac1{p^*}}\LEQ C\left(\sum_{i=1}^n\int_\O|\p_iv|^{p_i}\right)^{\frac1p}
\end{equation*}
for all $v\in W_0^{1,\vec p}(\O),\ if\ \DST\sum_{i=1}^n\dfrac1{p_i}>1.$
\item If $\DST\sum_{i=1}^n\dfrac1{p_i}<1$, then
$$W_0^{1,\vec p}(\O)\subset_{\!>} L^\infty(\O).$$
Moreover, there exists a constant $C(n)>0$ such that
$$||v||_\infty \LEQ C(n)\DST\prod_{i=1}^{n}\left\|\dfrac{\p v}{\p x_i}\right\|_{p_i}^{\frac {1}n}.$$
\item If $\DST\sum_{i=1}^n\dfrac1{p_i}=1$, then $$ W_0^{1,\vec p}(\O)\subset_{\!>} L^r(\O)$$ for all $r<+\infty.$
\end{enumerate}
\end{coro}
\begin{rem}
The two last statements can be found also  in \cite{Ta2}.
\end{rem}
As to the case of variable exponent spaces,   for $u:\O\to\R$ measurable, we  set
 $$\Phi_{p(\cdot)}(u)=\int_\O|u(x)|^{p(x)}dx$$
and we consider the norm:
\begin{equation}\label{eq1.4}
||u||_{p(\cdot)}=\inf\left\{\lambda>0:\Phi_{p(\cdot)}\left(\frac u\lambda\right)\LEQ1\right\},\qquad (\inf\emptyset=+\infty).
\end{equation}
Setting
$$L^{p(\cdot)}(\O)=\{u:\O\to\R\hbox{ measurable such that } ||u||_{p(\cdot)}<+\infty\},$$
the space $(L^\pp(\O);||\cdot||_\pp)$ is a Banach function space and an equivalent norm for $u$ is the following Amemiya norm
\begin{equation}\label{eq1.5}
|u|_\pp=\inf_{\lambda>0}\lambda\left(1+\Phi_{p(\cdot)}\left(\frac u\lambda\right)\right),
\end{equation}
which is equivalent to the norm in (\ref{eq1.4}) since
\begin{equation}\label{eq1.6}
||u||_\pp\LEQ|u|_\pp\LEQ2||u||_\pp.
\end{equation}
We set
$$ L^1_+(\O)=\{v\in L^1(\O): v\GEQ0\}\ ,\quad L^\pp_+(\O)=L^\pp(\O)\cap L^1_+(\O).$$
We {\bf  always assume} that
$$1<p_{-} = \inf \{ p(x): x\in \O\} \LEQ  p_{+} =\sup \{p(x): x\in \O\}\ < \infty.$$
\begin{propo}[\cite{CF}, Corollary 2.81 p. 63 and Corollary 2.23 p. 25]\label{p1112}\ \\
Under the above assumptions on $p$, one has:
\begin {itemize}
\item $L^{p(\cdot)}(\O)$ is reflexive.
 \item For all $u\in L^{\point}(\O)$,
$$||u||_\point\LEQ\left(\int_\O|u(x)|^{p(x)}dx\right)^{\frac1{p_-}}+\left(\int_\O|u(x)|^{p(x)}dx\right)^{\frac1{p_+}}.$$
\end{itemize}
\end{propo}
We also have a Poincar\'e-Sobolev type inequality for variable exponent spaces. Following
\cite{DHHR}, \cite{CF} for the next theorems (see also \cite{gdffio}),
we shall consider  exponents $p(\cdot)$   being bounded log-H\"older continuous functions on a bounded open set $\O$, i.e. satisfying the property\\
\it There exists a constant $c_1>0$ such that
 $$  \Log (e +1/|x-y| )|p(x)-p(y)| \LEQ c_1,
\forall (x,y) \in \O \times \O.$$
\rm Assuming also $p_+<n$, one can consider the Sobolev variable exponent
$$p^*(x)= \frac {np(x)}{n-p(x)},\;\; x \in \O,$$ and the following Poincar\'e-Sobolev
inequality holds true:
\begin{theo}\label{tVe1}\ \\
 There exists a constant $C=C(n,p(\cdot))$ such that
\begin{equation}\label{eqN300}
||v||_{p*(.)}\LEQ C || \nabla v||_{p(\cdot)}\qquad \hbox{for all $v\in W_0^{1,p(\cdot)}(\O)$.}
\end{equation}
\end{theo}
The dual of $W^{1,\point}_0$ is denoted by $W^{-1,p'(\cdot)}(\O)$. As usual, here $p'(x):=\dfrac{p(x)}{p(x)-1}$.

 We can summarize  the definitions of Lebesgue, Lorentz and Zygmund spaces as follows:
 \begin{defi}{\bf (Lorentz-Zygmund spaces) }\ \\
Let $\O$ be a space of measure  1, $0<p$, $q\LEQ +\infty,$ $-\infty<\lambda<+\infty.$ Then the Lorentz-Zygmund space $L^{p,q}\left(\log L\right)^{ \lambda }$ consists of all  Lebesgue measurable function $f$ on $\O$ such that :
$$||f||_{p,q;\lambda}=\begin{cases}\left(\DST\int_0^1\Big[t^{\frac 1 p}(1-\Log t)^{\lambda }f_*(t)\Big]^q\frac{dt}t\right)^{\frac1 q}&0<q<+\infty\\
\DST\sup_{0<t<1}t^{\frac 1p}(1-\Log t)^\lambda f_*(t)&q=+\infty\end{cases}\qquad\hbox{is finite.}$$
Here $f_*$ is the decreasing rearrangement of $|f|$.
\end{defi}
We also need the next definition.
\begin{defi}{\bf (of $G\Gamma(p,m;w_1,w_2)$) (see \cite{FFGKR})}\ \\
Let $w_1, w_2$ be two weights on $(0,1)$, $m\in[1,+\infty]$, $1\LEQ p<+\infty$. We assume the following conditions:
\begin{itemize}
\item[(c1)] There exists $K_{12}>0$ such that $w_2(2t)\LEQ K_{12}w_2(t)$ $\forall t\in(1/2,1)$. 
\item[(c2)] The function $\DST \int_0^tw_2(\s)d\s$ belongs to $L^{\frac mp}(0,1;w_1)$.
\end{itemize}
A generalized Gamma space with double weights is the set

$$G\Gamma(p,m;w_1,w_2)=\Big\{w:\O\to\R \hbox{ measurable and}\int_0^tv_*^p(\s)w_2(\s)d\s \hbox{ is in }L^{\frac mp}(0,1;w_1)\Big\},$$
which is a quasi-normed space endowed with the natural quasi-norm:
$$\rho(v)=\Big[ \int_0^1w_1(t)\Big(\int_0^tv_*^p(\s)w_2(\s)d\s \Big)^{\frac m p}dt \Big]^{\frac 1 p}.$$
If $w_2=1$ we simply denote $G\Gamma(p,m;w_1,1)=G\Gamma(p,m;w_1).$
\end{defi}
We shall also need the following elementary inequalities that can be found in \cite{LB},
\cite{DiazNPDE}.\\
For $p\GEQ2$,   there exists a constant $\a_p>0$ such that $\forall\,\xi\in\R^n,\ \forall\xi'\in\R^n$
\begin{equation}\label{eq1300}
\Big(|\xi|^{p-2}\xi-|\xi'|^{p-2}\xi',\xi-\xi'\Big)_{\R^n}\GEQ\a_p|\xi-\xi'|^p.\end{equation}
where in the left hand side the symbol $(\cdot,\cdot)_{\R^n}$ denotes the inner product in $\R^n$, and the symbol $|\cdot|$ is the associated norm.

A similar relation holds for the case $1<p <2$, namely, there exists
  a constant $\a_p>0$ such that $\forall\,\xi\in\R^n,\ \forall\xi'\in\R^n$
\begin{equation}\label{eq1301}
\Big(|\xi|^{p-2}\xi-|\xi'|^{p-2}\xi',\xi-\xi'\Big)_{\R^n}\GEQ\a_p\frac {|\xi-\xi'|^2}{(|\xi| +|\xi'|)^{p-2}}.\end{equation}

\section{\large \bf Abstract results on nonlinear interpolation}
We shall need the following results concerning real interpolation  with logarithm function (see \cite{EOP, GOT}).\\
Let $(X_0,||\cdot||_0),\ (X_1, ||\cdot||_1)$ be  two normed spaces  continuously embedded in a  Hausdorff topological vector space, that is, $(X_0,X_1)$ is  a compatible couple.
For $g\in X_0+X_1,\ t>0$, we shall denote
$$K(g,t)\dot=K(g,t;X_0,X_1)=\inf_{g=g_0+g_1}\Big(||g_0||_0+t||g_1||_1\Big).$$

For $0\LEQ\theta\LEQ1,\ \ 1\LEQ q\LEQ+\infty,\ \ \a\in\R,$ we define the interpolation space
$$(X_0,X_1)_{\theta,q;\a}=\Big\{ g\in X_0+X_1,\ ||g||_{\theta,q;\a}=||t^{-\theta-\frac1q}(1-\Log t)^\a K(\,g,t)||_{L^q(0,1)}\hbox{ is finite}\Big\}.$$

Next, we consider four normed spaces
$X_1\subset X_0,\, Y_1\subset Y_0$,   and ${\mathcal T}$ a non-linear mapping from $X_i$ into $Y_i$, $i=0,1$ such that:
\begin{enumerate}
\item $||{\mathcal T}a-{\mathcal T}b||_{Y_0}\LEQ f\Big(||a||_{X_0};||b||_{X_0}\Big)||a-b||^\a_{X_0}\hbox{ for all }(a,b)\hbox{ in }X_0.$
\item $||{\mathcal T}a||_{Y_1}\LEQ g\Big(||a||_{X_0}\Big)||a||_{X_1}^\beta,\ \ \forall\,a\in X_1.$
\end{enumerate}
Here $0<\a\LEQ 1$, $\beta>0$, $g$ is a continuous increasing function, and $f$ is continuous   on $\R^2$ and  such that for each $\sigma$,  $f(\sigma;\cdot)$ is increasing.\\

\subsection{\bf Estimating $K$-functional related to the mapping ${\mathcal T}$}
\begin{lem}\label{l1}\ \\
Under the above assumptions (1) and (2) on ${\mathcal T}$,
let $G(\sigma)=\Max\Big(g(2\sigma);f(\sigma;2\sigma)\Big),\sigma\in\R_+.$  Then for all $a\in X_0$, all $t>0$  one has
$$K\big({\mathcal T}a,t^\beta,Y_0,Y_1\big)=K({\mathcal T}a,t^\b)\LEQ G(||a||_{X_0})[K(a,t)^\beta+K(a,t)^\a].$$
Moreover, if $\beta\GEQ \alpha, $ then
$$K({\mathcal T}a,t^\b)\LEQ G(||a||_{X_0})(1+||a||_{X_0}^{\beta-\a})K(a,t)^\a. $$
\end{lem}
\Pr We follow Tartar's idea \cite{Tartar3} (see also \cite{Ma, Peetre}). If $a \in X_0$ and $ \eps > 0$, then there exist functions $a_0(\eps;\cdot)$
and $a_1(\eps,\cdot)$ such that
 $a_0(\eps,t)\dot=a_0(t)\in X_0,$ $a_1(\eps,t)\dot=a_1(t)\in X_1$ with
$a=a_0(t)+a_1(t)$ such that
\begin{equation}\label{eq1}
K(a,t)\LEQ||a_0(t)||_{X_0}+t||a_1(t)||_{X_1}\LEQ (1+\eps)K(a,t),\quad\forall\,t>0.
\end{equation}
We set ${\mathcal T}a=b_0(t)+b_1(t)$ with $b_1(t)={\mathcal T}a_1(t).$
Then
\begin{eqnarray}\label{eq2}
K({\mathcal T}a,t^\beta)&\LEQ &||b_0(t)||_{Y_0}+t^\b||b_1(t)||_{Y_1}=
||{\mathcal T}a-{\mathcal T}a_1(t)||_{Y_0}+t^\beta||{\mathcal T}a_1(t)||_{Y_1}
\nonumber\\&\LEQ& t^\b g\Big(||a_1(t)||_{X_0}\Big)||a_1(t)||_{X_1}^\b+f\Big(||a||_{X_0};||a_1(t)||_{X_0}\Big)||a_0(t)||^\a_{X_0}.
\end{eqnarray}
Since $a\in X_0, $ then
\begin{equation}
\label{eq3}
K(a,t)\LEQ ||a||_{X_0}.
\end{equation}
From relations (\ref{eq1}) and (\ref{eq3}), we have
\begin{equation}
\label{eq4}
||a_0(t)||_{X_0}\LEQ(1+\eps)||a||_{X_0}\quad\forall\,t>0,
\end{equation}and then
\begin{equation}\label{eq5}
||a_1(t)||_{X_0}\LEQ||a||_{X_0}+||a_0(t)||_{X_0}\LEQ(2+\eps)||a||_{X_0}.
\end{equation}
Therefore relation (\ref{eq2}) implies that
\begin{equation}\label{eq6}
K({\mathcal T}a,t^\beta)\LEQ\Max\Big(g\big(||a||_{X_0}(2+\eps)\big);
f\big(||a||_{X_0};(2+\eps)||a||_{X_0}\big)\Big)
\Big[||a_0(t)||^\a_{X_0}+t^\b||a_1(t)||^\b_{X_1}\Big],
\end{equation}
and combining this relation (\ref{eq6}) with relation (\ref{eq1}), and letting $\eps\to0$, we derive
\begin{equation}\label{eq7}
K({\mathcal T}a,t^\b)\LEQ G(||a||_{X_0})\Big[K(a,t)^\a+K(a,t)^\b\Big]\quad\forall\,t>0.
\end{equation}
If $\beta \GEQ \a$, then using relation (\ref{eq3}), one deduces from (\ref{eq7}) that
\begin{equation}\label{eq161}
K({\mathcal T}a,t^\b)\LEQ G(||a||_{X_0})(1+||a||_{X_0}^{\beta-\a})K(a,t)^\a.
\end{equation}
\ \HF
As  a particular case we have the following:
\begin{corol}[\bf of Lemma \ref{l1}]\label{c1l1} 
Let $X_1\subset X_0$, $Y_1\subset Y_0$ be four normed spaces. Assume that ${\mathcal T}:X_1\to Y_1$ is globally $\a$-H\"olderian, i.e. $\exists \, M_1>0$  such that
$$||{\mathcal T}a-{\mathcal T}b||_{Y_1}\LEQ M_1||a-b||^\a_{X_1},\quad 0<\a\LEQ1,$$
and ${\mathcal T}$ maps $X_0$ into $Y_0$ in the sense that $\exists\, M_0>0$, $\b >0$
$$||{\mathcal T}a||_{Y_0}\LEQ M_0||a||^\b_{X_0}.$$
Then, $\forall\,a\in X_0,\ \forall\,t>0$, one has
$$K({\mathcal T}a,t^\b)\LEQ\Max(M_0;M_1)\Big[K(a,t)^\b+K(a,t)^\a\Big].$$
If $\a\LEQ\b$, then
$$K({\mathcal T}a,t^\b)\LEQ (1+||a||_{X_0}^{\b-\a})\Max(M_0;M_1)K(a,t)^\a.$$
\end{corol}
\begin{corol}[\bf of Lemma \ref{l1}]\label{c2l1} 
Let $X_1\subset X_0, \ Y_1\subset Y_0$ be four normed  spaces. Assume that ${\mathcal T}:X_i\to Y_i$ is globally
$\a$-H\"olderian for $i=0,1$. Then, for all $a\in X_0, \ b\in X_1$ one has
$$K({\mathcal T}a-{\mathcal T}b;t^\a)\LEQ 2\Max(M_0;M_1) K(a-b;t)^\a,$$
where $M_i$ denotes the H\"older constants as in Corollary\ \ref{c1l1} of Lemma \ref{l1}.
Furthermore, if $X_1$ is dense in $X_0$, then the above equality holds also for
all $b \in X_0$.
\end{corol}
\Pr
Let $b\in X_1$ and define $Sx={\mathcal T}(b+x)-{\mathcal T}b$  for $x\in X_0$. Then
$$||Sx||_{Y_0}\LEQ M_0||x||^\a_{X_0},$$
and for all $x\in X_1$ and $y\in X_1$ we have
$$||Sx-Sy||_{Y_1}\LEQ M_1||x-y||^\a_{X_1}.$$
We may apply Corollary \ref{c1l1} of   Lemma \ref{l1} to derive
$$K(Sx;t^\a)\LEQ 2\Max(M_0;M_1)K(x;t)^\a, \forall\,x\in X_0.$$
Choosing    for $a\in X_0$, $x=a-b$ and taking into account that
$S(a-b)={\mathcal T}a-{\mathcal T}b$, we obtain  the first result.
If $X_1$ is dense in $X_0$, we consider a sequence $b_n \in X_1$ converging to $b$ in
$X_0$, and since  $$K(b_n -b;t) \LEQ ||b_n-b||_{X_0},$$ we have that $K(b_n -b;t)$ converges to zero as
$n$ goes to $\infty$. Writing
$$K({\mathcal T}a-{\mathcal T}b;t^\a)\LEQ K({\mathcal T}a-{\mathcal T}b_n;t^\a) + K({\mathcal T}b_n-{\mathcal T}b;t^\a),$$
and applying the preceding results, one has
$$K({\mathcal T}a-{\mathcal T}b_n;t^\a) + K({\mathcal T}b_n-{\mathcal T}b;t^\a) \LEQ 2\Max(M_0;M_1)\Big[K(a-b_n;t)^\a +K(b_n-b;t)^\a\Big].$$
Letting $n$ go to $\infty$, we get from the two last formulae that
$$K({\mathcal T}a-{\mathcal T}b;t^\a)\LEQ 2\Max(M_0;M_1) K(a-b;t)^\a, $$ for all $(a,b) \in X_0 \times X_0$.
\
\HF
\subsection{\bf Interpolation of H\"olderian mappings}
\begin{theo}\label{t1}\ \\
Let  $X_1\subset X_0$, $Y_1
\subset Y_0$, be four normed spaces, ${\mathcal T}$ be the mapping satisfying (1) and (2), and assume that $\a\LEQ\b$. Then, if $0\LEQ\theta\LEQ1,\ 1\LEQ p\LEQ+\infty$, for $a\in (X_0,X_1)_{\theta,p;\lambda}$ one has:
$${\mathcal T}a\in (Y_0,Y_1)_{\theta\frac\a\b,\frac p\a;\lambda\a}\qquad
\hbox{and}\qquad
||{\mathcal T}a||_{(Y_0,Y_1)_{\theta\frac\a\b,\frac p\a;\lambda \alpha}}\lesssim
\Big[(1+||a||_{X_0}^{\b-\a})G(||a||_{X_0})\Big]||a||^\a_{\theta,p;\lambda}.$$

\end{theo}
\Pr
One has from relation (\ref{eq161})
\begin{equation}\label{eq8}
K({\mathcal T}a,t^\b)\LEQ G(||a||_{X_0})(1+||a||_{X_0}^{\b-\a})K(a,t)^\a.
\end{equation}
Thus, by definition of $\|\cdot\|_{\theta,p;\lambda}$ (see the beginning of this section),
\begin{equation}\label{eq9}
J=\int_0^1t^{-\theta p}(1-\Log t)^{p\lambda}K({\mathcal T}a,t^\b)^{\frac p\a}\dfrac{dt}t\LEQ\Big[(1+||a||_{X_0}^{\b-\a})G(||a||_{X_0})\Big]^{\frac p\a}||a||^p_{\theta,p;\lambda}.
\end{equation}
Set $$J_1=\int_0^1\s^{-\theta\frac p\b}(1-\Log(\s))^{p\lambda}K({\mathcal T}a,\s)^{\frac p\a}\dfrac{d\s}\s .$$
We make the change of variables $\s=t^\b$ in the first integral $J $ to deduce:
\begin{equation}\label{eq10}
J=\frac1\beta\int_0^1\s^{-\theta\frac p\b}(1+\frac1\beta|\Log(\s)|)^{p\lambda}K({\mathcal T}a,\s)^{\frac p\a}\dfrac{d\s}\s.
\end{equation}
Hence $$ c_1J_1 \LEQ J \LEQ c_2J_1\, ,$$ \\with
$c_1=\begin{cases}\dfrac1\beta\Min\left(1;\dfrac1\beta\right)^{p\lambda}&\hbox{if $\lambda\GEQ 0,$}\\
\dfrac1\beta\Max\left(1;\dfrac1\beta\right)^{p\lambda}&\hbox{if $\lambda<0,$}\\
\end{cases}$ and $c_2=\begin{cases}\dfrac1\beta\Max\left(1;\dfrac1\beta\right)^{p\lambda}&\hbox{if $\lambda\GEQ 0,$}\\
\dfrac1\beta\Min\left(1;\dfrac1\beta\right)^{p\lambda}&\hbox{if $\lambda <0.$}\end{cases}$\\
\ \\We obtained that $J_1\simeq  J,$ and  from relations (\ref{eq9}), (\ref{eq10}) we conclude that
\begin{equation}\label{eq11}
||{\mathcal T}a||_{\theta\frac \a\b,\frac p\a;\lambda\a}\lesssim (1+||a||_{X_0}^{\b-\a})G(||a||_{X_0})||a||^\a_{\tpl}.
\end{equation}
\HF
In view of applications in P.D.E.s, we first have the following:
\begin{theo}\label{t2} \ \\
Let $X_1\subset X_0, \ Y_1\subset Y_0$ be four normed  spaces. Assume that ${\mathcal T}:X_i\to Y_i$ is globally
$\a$-H\"olderian for $i=0,1$.
For $0\LEQ\theta\LEQ1,\ 1\LEQ p\LEQ+\infty$, if $X_1$ is dense in $X_0$, then
${\mathcal T}$ is an $\a$-H\"olderian mapping from $(X_0,X_1)_{\theta,p;\lambda}$ into $(Y_0,Y_1)_{\theta,\frac p\a;\lambda \a}.$
\end{theo}
\Pr
Let $a \in (X_0,X_1)_{\theta,p;\lambda}$ and $b \in (X_0,X_1)_{\theta,p;\lambda}.$
 We have shown in the above Corollary \ref{c2l1} that $$K({\mathcal T}a-{\mathcal T}b;t^\a)\LEQ 2\Max(M_0;M_1) K(a-b;t)^\a. $$
 Following the same arguments as in proof of the above Theorem \ref{t1}, we deduce that
\begin{equation}\label{eq200}
||{\mathcal T}a-{\mathcal T}b||_\tpala\LEQ c_0||a-b||^\a_\tpl.\qquad\qquad
\end{equation}
\subsection{\bf Identifications of some interpolation spaces}
To obtain similar results as for Proposition \ref{p3} below with an interpolation process including a functor (as a logarithm function), we must identify the following interpolation spaces:
$$\Zue,\ \ \Zde,\ \ \Zte,\ \ \Zqe,$$
under suitable conditions on $s,\ p_2,\ \theta.$
Here is the general result collecting the necessary interpolation identification that we shall need for the application. The proof can be essentially found in (\cite{GOT}) (see also \cite{AFH}), however we give below the idea how to prove the statements.
\begin{propo}\label{p2.1}
Let $1\LEQ r<m\LEQ+\infty$,  $1\LEQ q_1$, $q_2\LEQ \infty$, $1\LEQ p<+\infty,$ $0\LEQ \theta<1$ and $\lambda\in\R$, if $\theta=1$ then $\lambda<-\dfrac 1p$, and $\lambda\GEQ - \dfrac 1p$ if $\theta=0$.
$$||f||_{(L^{r,q_1},L^{m,q_2})_{\theta,p,\lambda}}\simeq\begin{cases}
\Big[\DST\int_0^1\Big(t^{\frac {1-\theta}r+\frac\theta m}f_*(t)(1-\Log t)^ \lambda\Big)^p\dfrac {dt}t\Big]^{\frac1p}, &0<\theta<1;\\ \ \\
\Big[\DST\int_0^1\Big(\big(\int_0^t s^{\frac {q_1}r-1}f_*(s) ^{q_1}ds\big)^{\frac1{q_1}}(1-\Log t)^ \lambda\Big)^p\dfrac {dt}t\Big]^{\frac1p},
&\theta=0,\ q_1<+\infty;\\ \ \\
\Big[\DST\int_0^1\Big((\esssup_{0<s<t}s^{\frac 1r}f_*(s)\big) (1-\Log t)^ \lambda\Big)^p\dfrac {dt}t\Big]^{\frac1p},
&\theta=0,\ q_1=+\infty;\\ \ \\
\Big[\DST\int_0^1\Big(\big(\int_t^1 s^{\frac {q_2}m-1 }f_*(s)^{q_2}ds\big)^{\frac1{q_2}}(1-\Log t)^ \lambda\Big)^p\dfrac {dt}t\Big]^{\frac1p},
&\theta=1,\ q_2<+\infty;\\ \ \\
\Big[\DST\int_0^1\Big(\big(\esssup_{0<s<t}s^{\frac {1}m}f_*(s)) (1-\Log t)^ \lambda\Big)^p\dfrac {dt}t]\Big]^{\frac1p},
&\theta=1,\ q_2=+\infty.\\\end{cases}$$
\end{propo}
\Pr
The above proposition can be proved directly without invoking the general framework
in the cited references. Indeed, the main steps are the following two: first, to use Holmsted's formula
(see \cite{BenettSharpley}) to get an equivalent expression for the $K$-functional between Lorentz spaces and then to use suitable Hardy inequalities essentially developed in \cite{BR}.
\HF
We have several consequences of the above proposition. First, when we compare  the definitions of Lorentz-Zygmund spaces, $G\Gamma$-spaces, and small Lebesgue spaces (see \cite{FFGKR} and references therein for definition and properties), we have:
\begin{propo}\label{p2.2}\ \\
Let $1\LEQ r<m\LEQ\infty $, $\lambda\in\R$, $1\LEQ p_2<+\infty$.
\begin{enumerate}
\item If $0<\theta <1$, then the interpolation space $(L^{r;\infty};L^m)_{\theta,p_2;\lambda}$ coincides with the Lorentz-Zygmund space $L^{m_\theta,p_2}(\Log L)^\lambda$ with $\dfrac1{m_\theta}=\dfrac{1-\theta}r+\dfrac\theta m$.
$(L^1,L^m)_{\tpdl}$ coincides with the Lorentz-Zygmund space $L^{\frac{m'}{m'-\theta},p_2}\Big(\Log\,L\Big)^{\lambda}$.\\
\item If $\theta=0,\ 1\LEQ q_1,\ q_2\LEQ+\infty$, then for $j=1,2
$, we have
$$(L^{r,q_1},L^{m,q_2})_{0,p_2;\lambda}=G\Gamma(q_j,p_2;w_1,w_{2j}),$$
with $w_1(t)=(1-\Log(t))^{\lambda {p_2}}t^{-1}$, $w_{21}(t)=t^{\frac{q_1}r-1}$ if $q_1<\infty$, $w_{22}(t)=t^{\frac{q_2}m-1}$ if $q_2<\infty$, $t\in]0,1[$.
The space  $(L^1,L^m)_{0,p_2:\lambda}$ is the Generalized-Gamma space $G\Gamma(1,p_2;w_1)$ where $w(t)=t^{-1}(1-\Log t)^{\lambda p_2}$ (see \cite{FFGKR} or  \cite{FR}).
\item In particular, for $\lambda>-\dfrac1{p_2}$, we have the link with small Lebesgue spaces as follows:
\begin{enumerate}
\item If $1<q_1<+\infty,$ \qquad $L^{(q_1,\alpha_1}(\O)=(L^{q_1},L^{m,q_2})_{0,p_1;\lambda}$ \ \ $\forall\,q_2\in[1,+\infty],\ \forall\,m\in]q_1,+\infty],$\ \ $\alpha_1=q'_1(\lambda {p_2}+1)$,\ \
$\dfrac1{q_1}+\dfrac1{q'_1}=1.$\\
\item If $1<q_2<+\infty, $  \qquad $L^{( q_2,\alpha_2}(\O)=(L^{r,q_1},L^{q_2})_{0,p_2;\lambda}$ with $\alpha_2=q'_2(\lambda {p_2}+1), $ \\ $ 1=\dfrac1{q_2}+\dfrac1{q_2'},$ $ \forall\,r\in[1,q_2[,\ \ \forall\,q_1\in[1,+\infty]$.\\

\end{enumerate}
\end{enumerate}
\end{propo}

For the case $0<\theta<1$, we may apply the following duality result (see\cite{FFGKR, Peetre, BruK,Ta1}).
\begin{propo}\label{p6}\ \\
Let $X_1\subset X_0$  two Banach function spaces. Then the associate space of  $(X_0,X_1)_{\theta,p;\lambda}$ with $0<\theta<1,\ 1\LEQ p<+\infty$, $\lambda\in\R$, is the space
$$(X'_1,X'_0)_{1-\theta,p';-\lambda}\hbox{ with $\dfrac1p+\dfrac1{p'}=1$, where $X'_i$ is the associate space of $X_i,\ i=0,1.$}$$
\end{propo}
As a consequence we have the following
\begin{corop}[\bf of Proposition \ref{p6}]\label{c1p6}\ \\
Let $1<m<+\infty, \ 1<p'\LEQ+\infty,\ \lambda\in\R,\ 0<\theta<1$, $m'=\dfrac m{m-1}$. Then
$(L^m,L^\infty)_{\theta,p';\lambda}$ is the associate space of $(L^1,L^{m'})_{1-\theta,p';-\lambda},$
that is the Lorentz-Zygmund space (up to equivalence of norms) $(L^m,L^\infty)_{\theta,p';\lambda}=L^{\frac m{1-\theta},p'}\Big(\Log L\,\Big)^{ \lambda }.$
Moreover, we have
$(L^m,L^\infty)_{\theta,p';\lambda}=(L^{m,\infty},L^\infty)_{\theta,p';\lambda}.$
\end{corop}

Finally, we shall need the next result about a reiteration of Lorentz-Zygmund spaces, which follows from the Lions-Peetre's lemma (see  \cite{BenettSharpley, GOT}).
\begin{propo}\label{p10p10}(see \cite{GOT})\\
Let $1\LEQ p_0,\ q_0,\ p_1, q_1 \LEQ +\infty, p_1,\ 0<\theta<1, \ r_i\in \R$. Then
$$\left(L^{p_0,q_0}\Big(\Log L\Big)^{r_0},\ L^{p_1,q_1}\Big(\Log L\Big)^{r_1}\right)_{\theta,q;r}
=L^{p_\theta,q}\Big(\Log L\Big)^{r_\theta}$$
with $\dfrac 1{p_\theta}=\dfrac{1-\theta}{p_0}+\dfrac\theta{p_1}$\; and \;$r_\theta=(1-\theta)\dfrac{r_0 q}{q_0}+\dfrac{\theta r_1q}{q_1}.$
\end{propo}
Before starting the application of those interpolation formulae, we shall introduce  a very useful lemma inspired by the work of Benilan et al.  (see \cite{Benilanetal}, Lemma 4.2). But we state  it in a general framework in view of the applications to a large number of estimates that we shall use in the next section.
\section{\large \bf A fundamental lemma for estimates in Marcinkiewicz space}
\begin{lem}\label{lBt}{\bf Fundamental lemma of Benilan type}\\
Let $\nu $ be a non negative Borel measure and $h:\O\to\R_+,\ g:\O\to\R_+$, be two $\nu$-measurable functions. Then, $\forall\,\lambda>0,\ \forall\,k>0,$ we have\\
\C{$\DST\nu\Big\{h>\lambda\Big\}\LEQ\frac1\lambda \int_{\{g\LEQ k\}}hd\nu+\nu\{g>k\}.$}
\end{lem}
\Pr
Since $t\to\nu\{ h>t\}$ is non decreasing, $\forall\,\lambda>0$, we have
\begin{eqnarray*}
\nu\{h>\lambda\}&\LEQ&\frac1\lambda\int_0^\lambda\nu\{h>t\}dt\\
&=&\frac1\lambda\int_0^\lambda\Big(\nu\{h>t\}-\nu\{h>t,g>k\}\Big)dt
+\frac1\lambda\int_0^\lambda\nu\{f>t,g>k\}dt.
\end{eqnarray*}
We have $\nu\{h>t,g>k\}\LEQ\nu\{g>k\}$ and \\
\C{$\nu\{h>t\}-\nu\{h>t:g>k\}=\nu\{h>t:g\LEQ k\}.$}\\
Therefore, we obtain\\
\C{$\DST\nu\{h>\lambda\}\LEQ\frac1\lambda\int_0^{+\infty}\nu\{h>t:g\LEQ k\}dt+\nu\{g>k\}.$}\\
By the Cavalieri's principle, one has
$$\int_0^{+\infty}\nu\{h>t:g\LEQ k\}dt=\int_{\{g\LEQ k\}}hd\nu.$$
With those two last inequalities, we get the result.
\HF
Besides the applications that we shall give in the next section, we recall some estimates that we have already used in a previous work (\cite{DGRT}).
Let us recall that if $\omega$ is an integrable weight function on $\O$, the weighted Marcinkiewicz space is defined by
$$L^{q,\infty}(\O,\omega)=\left\{ v:\O\to\R\hbox{ measurable s.t. } \sup_{\lambda>0}\lambda^q\int_{|v|>\lambda}\omega(x)dx<+\infty\right\}, 0<q<+\infty.$$
If $\omega=1$, $L^{q,\infty}(\O,1)$ is the same as the Lorentz space $L^{q,\infty}(\O)$ defined in the first section.
\begin{theo}\label{tGB}{\bf A generalized Benilan type result}\\
Let $\omega$ be an integrable weight function on $\O$, $1\LEQ p<+\infty$, and let $u\in W^{1,1}_{loc}(\O)$ be such that, for a constant $M>0$, we have

$$\int_\O|\nabla T_k(u)|^p\omega(x)dx +\int_\O|T_k(u)|^p\omega(x)dx\LEQ Mk,\ \forall\,k>0,$$
with $T_k(t)=\Min(|t|;k) \sign(t),$  $t\in\R$. \\
Assume furthermore that we have a continuous Sobolev embedding,
$$W^1L^p(\O,\omega)\subset_{\!\!>} L^{p^*}(\O,\omega)\quad\hbox{ for some } p^*>p.$$
Then, one has:
$$\int_{|\nabla u|>\lambda}\omega dx\LEQ cM^{\frac{p^*}{p^*+p'}}\lambda^{-\frac{pp^*}{p^*+p'}}\quad\forall\,\lambda>0, $$
\hbox{ where $c>0$ is a constant depending only on $p,\ \O,\ p^*$, $p'.$}
Hence $|\nabla u|\in L^{q,\infty}(\O,\omega)$, with $q=\dfrac{pp^*}{p^*+p'}$.
If $q>1, $ then $u\in W^1L^r(\O,\omega),\ 1\LEQ r<q.$
\end{theo}
\Pr
For a measurable set $E\subset\O$, we set $\nu E=\DST\int_E\omega(x)dx$ and we apply the above fundamental Lemma \ref{lBt} to derive that for all $\lambda>0$, $\forall\,k>0$,
\begin{equation}\label{eq3901}
\nu\{|\nabla u|^p>\lambda\}\LEQ \dfrac1\lambda \int_\O|\nabla T_k(u)|^p\omega dx+\nu\{|u|>k|\}.
\end{equation}
By the first assumption of the  theorem, we get, for all $k>0$ and $\lambda>0,$ that
\begin{equation}\label{eq3902}
\nu\{|\nabla u|^{p}>\lambda\}\LEQ\dfrac k\lambda M+\nu\{|u|>k\}.
\end{equation}
Following Benilan et al. \cite{Benilanetal}, we have $\{|u|>\eps\}=\big\{|T_k(u)|>\eps\big\}$ for $\eps<k$. Therefore
\begin{equation}\label{eq3903}
\nu\{|u|>\eps\}\LEQ\dfrac1{\eps^{p^*}}\int_\O|T_k(u)|^{p^*}\omega dx.
\end{equation}
Using  Sobolev's inequality, we have
\begin{equation}\label{eq3904}
\nu\{|u|>\eps\}\LEQ\dfrac1{\eps^{p^*}}c_s\left[\int_\O|\nabla T_k(u)|^p\omega dx+\int_\O|T_k(u)|^p\omega dx\right]^{\frac{p^*}p}\LEQ\dfrac1{\eps^{p^*}}c_s k^{\frac{p^*}p}.
\end{equation}
As $\eps\to k$, we have, for all $k>0,$
\begin{equation}\label{eq3905}
\nu\{|u|>k\}=\int_{|u|>k}\omega(x)dx\LEQ c_sk^{-\frac{p^*}{p'}},
\end{equation}
where $c_s$ is the Sobolev's constant.
Combining relations (\ref{eq3902}) and (\ref{eq3905}), one has
\begin{equation}\label{eq3906}
\nu\{|\nabla u|^p>\lambda\}\LEQ \Inf_{k>0}\left\{\dfrac{M}\lambda k+c_sk^{-\frac{p*}{p'}}\right\}.
\end{equation}
Computing the infimum, we have $\forall\,\lambda>0$
\begin{equation}\label{eq3907}
\nu\{|\nabla u|^p>\lambda\}\lesssim M^{\frac{p^*}{p^*+p'}}\lambda^{-\frac{p^*}{p^*+p'}},
\end{equation}
and this implies the result.\HF
Here are some weighted spaces in which we have a Sobolev embedding (see \cite{Kufner, Simon}).
\begin{propo}\label{p3900}\ \\
Assume that $\O$ is a bounded open Lipschitz set of $\R^n$. Let $\a\GEQ0$,and let $\omega$ be one of the following weights
\begin{itemize}
\item $\omega(x)=\dist(x;\p\O)^\alpha=\delta(x)^\alpha$,
\item $\omega(x)=\dist(x;x_0)^\alpha,\quad x_0\in\p\O$.
\end{itemize}
For $1\LEQ p<n+\alpha$, we have $p^*=\dfrac{(n+\alpha)p}{n+\alpha-p}$ and
$$\left[\int_\O|v|^{p^*}\omega dx\right]^{\frac1{p^*}}\LEQ c\left[\left(\int_\O|v|^p\omega dx\right)^{\frac1p}+\left(\int_\O|\nabla v|^p\omega dx\right)^{\frac1p}\right].$$
\end{propo}
As a consequence of the above Theorem and Proposition \ref{p3900}, here is a proposition that we have already stated and used in \cite{DGRT} (see Proposition 13 therein).

\begin{propo}\label{p11}\ \\
	Let $v\in L^1(\O,\delta^\a),$ and $ \a\in[\,0,1\,]$. Assume that there exists a constant $c_0>0$ such that  for all $k>0$
	\begin{equation*}T_k(v):=\Min(|v|;k)\ {{\rm sign\,}}(v)\in W^1L^2(\O,\delta^\a),\end{equation*}
	and
	\begin{equation}\label{eq24bis}
	\int_\O|\nabla T_k(v)|^2\delta^\a dx+\int_\O|T_k(v)|^2\delta^\a dx\LEQ c_0 k.
	\end{equation}
	Then there exists a constant $c$, depending continuously on $c_0>0$, such that for all $\lambda>0$
	\begin{equation*}\int_{\{x:|\nabla v|(x)>\lambda\}}\delta^\a(x)dx\LEQ\dfrac c{\lambda^{1+\frac1{n+\a-1}}}.\end{equation*}
	
	In particular, if $(v_j)$ is a sequence converging weakly in $L^1(\O)$ to a function $v$, satisfying the inequality (\ref{eq24bis}) and such that
	\begin{equation*}\int_\O|\nabla T_k(v_j)|^2\delta^\a dx\LEQ c_0k\qquad\qquad\forall j,\ \forall k,\end{equation*}
	then $(v_j)$ converges to $v$ weakly in $W^{1,q}(\O')$ for all  $q\in\left[1,\dfrac{n+\alpha}{n+\a-1}\right[$ and all $\O'\subset\!\subset\O$, and there exists a subsequence (that we call still $(v_j)$) such that  $v_j(x)\to v(x)$ a.e. in $\O$.
\end{propo}
\section{\bf \large Application to the regularity of the solution of a  $p$-Laplacian}
Let $\O$ be a bounded set of $\R^n$. Let us consider $f\in L^1(\O)\cap W^{-1,p'}(\O),\ \dfrac 1p+\dfrac 1{p'}=1,\ 1<p<+\infty$, and
 $V$ a Caratheodory function from $\O\times \R$ into $\R$ such that
 \begin{description}
\item[(H1)] for all $\s\in\R, \ x\in\O\to V(x;\s)$ is in $L^\infty(\O)$.
\item[(H2)] for a.e. $x\in \O, \s\in\R\to V(x;\s)$ is continuous and non decreasing with $V(x;0)=0.$
\end{description}
 Using the Leray-Lions' method for monotone operators (see Lions's book \cite{Lionslivre}) or the usual fixed point theorem of Leray-Schauder's type (see Gilbarg -Trudinger \cite{GT})
  we have:
  \begin{propo}\label{p3001-s}\ \\
  Let $f$  be in $L^1(\O)\cap W^{-1,p'}(\O).$ Then there exists  a unique element
  $u\in W^{1,p}_0(\O)$ such that
\begin{equation}\label{eq12}
\int_\O\f(x)V(x;u)dx+\int_\O|\nabla u|^{p-2}\nabla u\cdot \nabla \f dx=\int_\O f\f\,dx\quad\forall\,\f\in W^{1,p}_0(\O)\cap L^\infty(\O).
\end{equation}
We call such solution a weak solution of the Dirichlet equation $-\Delta_pu+V(x;u)=f$.
\end{propo}
\begin{rem}{(\bf on the above existence and uniqueness)}\\
If $p>n$, then $L^1(\O) \subset W^{-1,p'}(\O)$.
If $p\LEQ n$, then the dual space $L^{(p^*)'}(\O)\subset W^{-1,p'}(\O)$ whenever $p^*=\dfrac{np}{n-p} $ if $p<n, $ and $p^*$ is any finite number if $p=n.$\\
Therefore, the above result can be applied for these cases.
In the paragraph concerning the equation with variable exponents, we give the idea on how
to prove the above proposition.
\end{rem}
We can define a nonlinear mapping:
$$\begin{matrix}
{\mathcal T}:&L^1(\O)\cap W^{-1,p'}(\O) &\longrightarrow&[L^p(\O)]^n\\
&f&\longmapsto&{\mathcal T}f=\nabla u.\end{matrix}$$
We shall need sometimes the following  additional growth assumption for $V.$
\begin{description}
\item[(H3)] There exist $m_1\in \Big[p-1,\OV m_1\Big[$, $\OV m_1=\begin{cases}(p-1)\left(1+\dfrac1{n-p}\right)&if\ p<n\\
<+\infty&if\ p\GEQ n,\end{cases}$\\and a constant $c>0$ such that
$$|V(x,\s)|\LEQ c|\s|^{m_1}\qquad\forall\,\,\s\in\R,\ a.e.\ x\in\O.
$$
\end{description}
We want to extend the above mapping over all $L^1(\O)$.
When $p=n$ and $ f\in L^1(\O)$, the Iwaniec-Sbordone's method ensures the existence and uniqueness of a weak solution that is  under the above formulation (\ref{eq12})
or even in the sense of distribution, see for instance (\cite{FS},   \cite{MPR},   \cite{FJRV}).
So the above mapping is well defined on $L^1(\O)$.\\
When $p<n$ and the data $f$ is only in $L^1(\O)$, the formulation by equation (\ref{eq12}) cannot ensure the uniqueness of the solution. Here it is an equivalent formulation which summarizes various definitions introduced by  different authors (see for instance
\cite{DipernaLions,Benilanetal,BlaMu,BGDM,Rako5,Rako7,Carilloetal}).
We consider again the usual truncation
$$T_k:\R\to\R\hbox{ defined by }T_k(\s)=\{|k+\s|-|k-\s|\}/2,$$
and we define as in \cite{Rako3,Rako5} (see also \cite{Benilanetal}), the following $T$-space or $T$-set:
$$\S^{1,p}_0=\Big\{v:\O\to\R\hbox{ measurable such that  $tan^{-1}(v)\in W^{1,1}_0(\O)$},$$
$$\hbox{ and  $\forall\,k>0,$  $T_k(v)\in W^{1,p}_0(\O)$, }\      \sup_{k>0}k^{-\frac1p}||\nabla T_k(v)||_{L^p(\O)}=\kappa<+\infty\Big\}.$$
\begin{defi}[\bf of an entropic-renormalized solution]\label{d1111}\ \\
We will say that a function $u$  defined  on $\O$ is an entropic-renormalized solution associated to the Dirichlet problem
\begin{equation}\label{eq14}\ \\
-\Delta_pu+V(x;u)
=f\in L^1(\O)\quad\hbox{ 	$u=0$ on $\p \O$ }
\end{equation}
if
\begin{enumerate}
\item $u\in\S^{1,p}_0(\O)$, $V(\cdot, u)\in L^1(\O).$
\item $\forall\,\eta\in W^{1,r}(\O),\ r>n,\ \ \forall\,\f\in W^{1,p}_0(\O)\cap L^{\infty}(\O)$ and all $B\in W^{1,\infty}(\R)$ with $B(0)=0$, $B'(\s)=0$ for all $\s$ such $|\s|\GEQ \s_0>0$, one has:
\begin{equation}\label{eq15}
\int_\O|\nabla u|^{p-2}\Nu\cdot\nabla\Big(\eta B(u-\f)\Big)dx+
\int_\O V(x;u)\eta  B(u-\f)dx=\int_\O f\eta B(u-\f)dx.
\end{equation}
\end{enumerate}
\end{defi}
If $f\in L^{p'}(\O)$ the  above formulation (\ref{eq15}) is equivalent to the formulation (\ref{eq12}) (i.e. a weak solution is an entropic-renormalized solution); see \cite{Rako3}
or \cite{MPR} for the case $p=n$.
It has been proved in the above references (see \cite{Benilanetal,Rako6,Rako7,BlaMu,Carilloetal} ) that we have existence and uniqueness of an entropic-renormalized solution.
\begin{theo}\label{t3.1}\ \\
Let $f\in L^1(\O)$ and assume (H1) and (H2). Then there exists a unique entropic-renormalized solution of equation (\ref{eq14}). Moreover, if the sequence $(f_j)$ converges to $f$ in $L^1(\O)$,  then the sequence $(\Nu_j(x))$ converges to $\Nu(x)$ almost everywhere in $\O$ for a subsequence still denoted by $(u_j)$.
When $p>2-\dfrac 1n$, the solution $u\in  W^{1,1}_0(\O).$\\

\end{theo}
{\bf Comments on the proofs of Theorem \ref{t3.1} and Theorem \ref{t3.21}}\\
\begin{itemize}
\item For any $v\in\S^{1,p}_0(\O)$, the gradient of $v$ exists a.e. in the sense that if we denote by $\{\vec e_1,\ldots,\vec e_n\Big\}$ the canonical basis of $\R^n$, then the following limit exists
 almost everywhere in $\O$ $$\lim_{t\to0}\dfrac{v(x+t\vec e_i )-v(x)}t\dot=\dfrac{\p v}{\p x_i}(x)$$
and
$$DB(v)(x)=B'\big(v(x)\big)Dv(x),\hbox{ whenever }B\in W^{1,\infty}(\R).$$
This result is only given in \cite{Rako7} (see also \cite{Rako5}).
\item Let $f_1\in L^1(\O),\ f_2\in L^\infty(\O),\ u_1$ be the entropic-renormalized solution associated to $f_1$ and $u_2$ be the weak solution of equation (\ref{eq12}) associated to $f_2$. Then, choosing $\eta=1$,  $B=tan^{-1}(T_k),$ for $k>0$, $\f=u_2$ is in $W^{1,p}_0(\O)\cap L^\infty(\O). $ One has, dropping the non negative term,
\begin{equation}\label{eq16}
\int_\O\Big(|\nabla u_1|^{p-2}\nabla u_1-|\nabla u_2|^{p-2}\nabla u_2\Big)\cdot\nabla B(u_1-u_2)dx\LEQ\dfrac\pi2\int_\O|f_1-f_2|dx.
\end{equation}
The relation (\ref{eq16}) implies the uniqueness of the entropic-renormalized solution for all $p\in]1,n[$.\\
When $p\GEQ2$, we can have more inequalities for $u_1- u_2$. Indeed, we can
use the strong coercivity of the $p$-Laplacian, that is inequality
(\ref{eq1300}) (or see  below (\ref{eq13})), and we let $k\to+\infty$ to obtain:
$$\int_\O\dfrac{|\nabla(u_1-u_2)|^p}{1+|u_1-u_2|^2}dx\LEQ\dfrac\pi2\int_\O|f_1-f_2|dx.$$
From this relation, we have for all $1\LEQ q<\dfrac n{n-1}(p-1)$
$$\int_\O|\nabla(u_1-u_2)|^q\LEQ c||f_1-f_2||^{\frac q p}_{L^1}\left[1+\left(\int_\O|\nabla(u_1-u_2)|^qdx\right)^{m_2}\right]
$$ with $m_2=q^*\left(\dfrac1q-\dfrac1p\right)$ and $\dfrac1{q^*}=\dfrac1q-\dfrac1n$  if $p<n$,
so that $m_2<1$; and $q^*$ is any number so that $m_2<1$ if $p\GEQ n.$\\
Hence, using Young's inequality,
\begin{equation}\label{eq17}
\left(\int_\O|\nabla(u_1-u_2)|^q\right)^{1/q}\LEQ
c\Big[||f_1-f_2||^{\beta_1}_{L^1}+||f_1-f_2||^{\b_2}_{L^1}\Big]
\end{equation}
where $\b_i,\ i=1,2,\ c$ depend only $p,\ n,\ \O$, $\b_i>0$.
This is the method used in \cite{Rako2,Rako4}. The above inequality gives    a stability and uniqueness result.\\
The technique developed by Benilan et al. \cite{Benilanetal} gives a more precise
result than the above relation (\ref{eq17}). Indeed the same arguments as for having (\ref{eq16}) with $B=T_k$, leads to:
\begin{equation}\label{eq18}
\int_\O|\nabla T_k(u_1-u_2)|^pdx\LEQ k\int_\O|f_1-f_2|dx,\qquad\forall\, k\GEQ0.
\end{equation}
If $f_1\neq f_2$, we set $w=\dfrac{u_1-u_2}{||f_1-f_2||^{\frac1{p-1}}_{L^1}}.$ \ \
$\lambda=k||f_1-f_2||^{-\frac1{p-1}}_{L^1}$ and we deduce that
\begin{equation}\label{eq19}
\int_\O|\nabla T_\lambda(w)|^pdx\LEQ\lambda.
\end{equation}
From this inequality, Benilan's technique (see \cite{Benilanetal}  or  the above Theorem \ref{tGB}) implies that
\begin{equation}\label{eq20}
||\nabla  w||_{L^{n'(p-1),\infty}}\LEQ c(p,\O).
\end{equation}
This implies the second statement of the Theorem \ref{t3.21}. Another proof of this regularity result (\ref{eq20}) is in \cite{RakoHal}.\ \HF
\end{itemize}

\begin{propo}\label{p1}\ \\
Let $u$ be the solution of equation (\ref{eq14}) with $f$ being in $L^{p'}(\O)$. 
\begin{enumerate}
\item If $f \in L^{\frac np;\frac1{p-1}}(\O), p\LEQ n $, \hbox{then} $u\in L^\infty(\O)$ and $||u||_{L^\infty}\LEQ c||f||^{\frac1{p-1}}
 _{L^{\frac np;\frac1{p-1}}(\O)}\,$,\\
 if $f \in L^{1;\frac1{p-1}}(\O), p> n $, \hbox{then} $u\in L^\infty(\O)$ and $||u||_\infty \LEQ c||f||_{L^{1;\frac1{p-1}}(\O)}^{\frac 1{p-1}}$. 

\item If we assume (H3)  and $f \in L^{n,1}(\O)$, then
$$\nabla u\in L^\infty\hbox{ and }
||\nabla u||_{L^\infty}\LEQ c\Big(1+||f||^{\frac{m_1+1-p}{(p-1)^2}}_{L^1}\Big)||f||^{\frac1{p-1}}_{L^{n,1}(\O)}.$$
\end{enumerate}
All the constants denoted by $c$ depend only on $p$, $\O$ and $V$.
\end{propo}
This proposition gathers well known results (see for instance \cite{Cianchi-Mazja}, \cite{FPR} or \cite{Rakobook} page 125  for statement (1), and \cite{Cianchi-Mazja, Cianchi-Mazja2, RakoHal} for statement (2)). 
The growth of the gradient in Proposition \ref{p1} comes from the following Lemmas.
\begin{lem}\label{l3.1}\ \\
Assume (H3).
Let $m_3=\dfrac n{n-p}(p-1)$ if $1<p<n$, or $m_3$ be any finite number in $[nm_1,+\infty[$ if $p\ge n$, and let $u\in L^\infty(\O)$, then
$$||V(\cdot;u(\cdot))||_{L^{n,1}}^{\frac1{p-1}}\lesssim||u||_\infty\cdot||u||_{L^{m_3,1}}^{\frac{m_1+1-p}{p-1}}.$$
\end{lem}
\Pr
One has from (H3)
$$||V(\cdot;u(\cdot))||_{L^{n,1}}^{\frac1{p-1}}\LEQ c\left[\DST\int_0^{|\O|}t^{\frac1n}|u|_*^{m_1}(t)\dfrac{dt}t\right]^{\frac1{p-1}}$$
$$\ \ \lesssim||u||_\infty\left[\int_0^{|\O|}t^{\frac1n}|u|_*^{m_1+1-p}(t)\dfrac{dt}t\right]^{\frac1{p-1}}$$
\begin{equation}\label{eq310}
\lesssim||u||_\infty||u||_{L^{n(m_1+1-p),m_1+1-p}}^{\frac{m_1+1-p}{p-1}}.
\end{equation}
If $1<p<n$, one has $n(m_1+1-p)<\dfrac n{n-p}(p-1)=m_3$ since $m_1<(p-1)\left[1+\dfrac1{n-p}\right]$, \ $n(m_1+1-p)<m_3$ if $p\GEQ n$. Therefore, the last inequality implies the result.
\HF
\begin{lem}\label{l3121}\ \\
Assume  (H3) and  let $m_3$ be as in Lemma \ref{l3.1}. If $f\in L^{n,1}(\O)$ and $u$ is  a weak solution of equation (\ref{eq12}), then
$$||V(\cdot;u(\cdot))||^{\frac1{p-1}}_{L^{n,1}}\lesssim ||f||_{L^{n,1}}^\frac1{p-1}\cdot||f||_{L^1}^{\frac{m_1+1-p}{p-1}}.$$
\end{lem}
\Pr
By statement (1) of Proposition \ref{p1}, $u$ is bounded and $||u||_\infty\lesssim ||f||^{\frac1{p-1}}_{L^{n,1}}$. We have
\begin{equation}\label{eq312}
||u||_{L^{n'(p-1)^*,1}}\LEQ c||\nabla u||_{L^{n'(p-1),\infty}}
\end{equation}
by Sobolev-Poincar\'e's inequality (see \cite{Rakobook}), and
\begin{equation}\label{eq313}
||\nabla u||_{L^{n'(p-1),\infty}}\lesssim ||f||^{\frac1{p-1}}_{L^1}
\end{equation}
by Theorem \ref{t3.1}.
From (\ref{eq312}) and (\ref{eq313}), we have
\begin{equation}\label{eq314}
||u||_{L^{n'(p-1)^*,1}}\lesssim ||f||^{\frac1{p-1}}_{L^1}.
\end{equation}
Since  $n'(p-1)^*=m_3$,   from Lemma \ref{l3.1} we derive the result.
\HF
Statement (2) of Proposition \ref{p1} is then a consequence of the above Lemmas since  $-\Delta_p u=f-V(\cdot;u)\in L^{n,1}$ if $u$ is a weak solution of (\ref{eq12}) for $f\in L^{n,1}(\O).$ The Cianchi-Maz'ja's result implies
\begin{equation}\label{eq325}
||\nabla u||_{L^\infty}\lesssim||f-V(\cdot;u)||^{\frac1{p-1}}_{L^{n,1}}\lesssim||f||^{\frac1{p-1}}_{L^{n,1}}+||V(\cdot;u)||_{L^{n,1}}^{\frac1{p-1}}.
\end{equation}
This inequality (\ref{eq325}) and Lemma \ref{l3121} implies the estimate (2) in Proposition \ref{p1}.\\
\subsection{\bf  The H\"olderian mappings for the case $p\GEQ 2$}\ \\
We start with the H\"older property in the {\bf case $2\LEQ p < n$}.
\begin{theo}\label{t3.21}\ \\
 If $p\GEQ2$, $f_i\in L^1(\O),\ i=1,2$ and if $u_i,\ i=1,2$ are the corresponding entropic-renormalized solution, \begin{enumerate}
\item $\DST \int_\O\Big|\nabla T_k(u_1-u_2)\Big|^{p}dx\LEQ k\int_\O|f_1-f_2|dx.$
\item $u_i\in W_0^{1,1}(\O)$, and  moreover
$$||\nabla(u_1-u_2)||_{L^{n'(p-1),\infty}(\O)}\LEQ c||f_1-f_2||^{\frac1{p-1}}_{L^1(\O)},$$
where $c$ is a constant depending only on the data $p,\ \O$ and $V,\ n'=\dfrac n{n-1}.$
\end{enumerate}
\end{theo}
\Pr
We use the stability result. Indeed, let $f_{1j}=T_j(f_1)$ (resp. $f_{2j}=T_j(f_2)$). Then, $u_{i_j}$ $i=1,2$ associated to $f_{ij}$ are solutions of (\ref{eq12}).
We note first that if $u_1$ is an entropic-renormalized solution associated to $f_1$ and $f_{1j}=T_j(f_1)\in L^\infty(\O),$ then the weak solution $u_{1j}$ of equation (\ref{eq12}) satisfies
$$||\nabla u_{1j}-\nabla u_1||_{L^{n'(p-1),\infty}(\O)}\xrightarrow[j\to+\infty]{}0.$$
 Therefore  using relation (\ref{eq1300}) and $\f=T_k(u_{1{_j}}-u_{2_{j}})$ as a test function in relation (\ref{eq12}), we derive the statement (1) of the theorem using
the convergences for each $\Nu_i$.
While for statement (2), we may apply Theorem \ref{tGB} with $u=u_1-u_2\in W^{1,1}_{loc}(\O)$ and $\omega=1.$
\HF

As consequence of this theorem, we have the following Corollary which proves Theorem \ref{t3.1}.
\begin{coro}[\bf of Theorem \ref{t3.21}]\ \\
Under the assumptions (H1) and (H2), we extend the mapping ${\mathcal T}:\begin{matrix}L^1(\O)&\longrightarrow&\Big[L^{n'(p-1),\infty}(\O)\Big]^n\\
f&\longmapsto&{\mathcal T}f\end{matrix}$ with ${\mathcal T}f=\nabla u,$ where  $\ u$ is the unique entropic-renormalized solution of the Dirichlet equation (\ref{eq14}). Then,  for $p\GEQ2$, there exists a constant $c(p,\O)>0$ independent of $V$ such that
$$||{\mathcal T}f_1-{\mathcal T}f_2||_{L^{n'(p-1),\infty}}\LEQ c(p,\O)||f_1-f_2||^{\frac1{p-1}}_{L^1(\O)}.$$
\end{coro}

We derive this result from the  statement (2) of Theorem \ref{t3.21}. This stability implies the desired result.

\begin{lem}\label{l2}\ \\
Assume (H1) and (H2).
If $p\GEQ2$, then the preceding mapping ${\mathcal T}$ is $\dfrac1{p-1}$-H\"olderian  from $L^{(p^*)'}$ into $[L^p(\O)]^n$ with
$\dfrac1{(p^*)'}=\dfrac1{p'}+\dfrac1n,\quad p^*=\dfrac{np}{n-p},\quad 1<p<n$, $p'$ is the conjugate of $p$.
\end{lem}
\Pr
For $p\GEQ2$, we recall that there exists a constant $\a_p>0$ such that $\forall\,\xi\in\R^n,\ \forall\xi'\in\R^n$
\begin{equation}\label{eq13}
\Big(|\xi|^{p-2}\xi-|\xi'|^{p-2}\xi',\xi-\xi'\Big)_{\R^n}\GEQ\a_p|\xi-\xi'|^p.\end{equation}
Therefore, for two data $f_1$ and $f_2$ in $L^{p'}(\O)$, dropping the non negative term, we have
$$c_p\int_\O|{\mathcal T}f_1-{\mathcal T}f_2|^pdx\LEQ\int_\O\Big(|\nabla u_1|^{p-2}\nabla u_1-|\nabla u_2|^{p-2}\nabla u_2,\nabla(u_1-u_2)\Big)dx$$
$$\LEQ \int_\O(f_1-f_2)(u_1-u_2)dx\quad\hbox{(by Poincar\'e-Sobolev inequality)}$$
$$\LEQ c_{1p}||f_1-f_2||_{L^{(p^*)'}}||{\mathcal T}f_1-{\mathcal T}f_2||_{L^p},$$
so that
$$||{\mathcal T}f_1-{\mathcal T}f_2||_{L^p(\O)}\LEQ c_{2p}||f_1-f_2||^{\frac1{p-1}}_{L^{(p^*)'}(\O)}.$$
This implies the result using a density argument.
\HF
To apply the abstract results given in the second section for interpolation spaces, we need to use some well-known results concerning some identification. The first one can be deduced from the famous reiteration process of Lions-Peetre or from Proposition \ref{p2.1}.
\begin{propo}\label{p2}\ \\
For all $r\in [1,+\infty],$   $1\LEQ m\LEQ +\infty,$
$1<q\LEQ+\infty,$ $m<k<q$, we have
$$\Big(L^m(\O),\ L^q(\O)\Big)_{\theta,r;0}=L^{k,r}(\O)\hbox{ with }\dfrac1k=\dfrac{1-\theta}m+\dfrac \theta q.$$
\end{propo}
Notice that the interpolation space $(X_0,X_1)_{\theta,r;0}$ is the same as Peetre interpolation space $(X_0,X_1)_{\theta,r}$ since $X_1\subset X_0.$
\begin{propo}\label{p3}\ \\
Assume (H1) and (H2).
Let $p^*=\dfrac {np}{n-p}$ with $2\LEQ p<n$, $(p^*)'$ its conjugate,\\
 $1\LEQ k\LEQ (p^*)'=\dfrac{np}{np-n+p},$\ \  $p^*$ is any finite number if $p\GEQ n$, $r\in[1,+\infty]$. Then
$$||{\mathcal T}f_1-{\mathcal T}f_2||_{L^{k^*(p-1),r(p-1)}}\LEQ c||f_1-f_2||_{L^{k,r}}^{\frac1{p-1}},$$
for $f_1,\ f_2$ in $L^{k,r}(\O)$, with $k^*=\dfrac{kn}{n-k}$ if $k<n$ and  any finite number if $k\GEQ n.$\\
In particular if $f\in L^{k,r}(\O)$ then the gradient of the solution $u$ of  (\ref{eq15}) belongs to $[L^{k^*(p-1),r(p-1)}(\O)]^n$.
\end{propo}
\Pr
The mapping ${\mathcal T}$ is $\dfrac1{p-1}$-H\"olderian from $L^1(\O)$ into $[L^{n'(p-1),\infty}(\O)]^n$ and from $L^{q'}(\O)$ into $[L^p(\O)]^n$ with $q'=(p^*)'$. Moreover,  we have  $L^{k,r}(\O)=(L^1,L^{q'})_{\theta,r}$ with $\theta=p^*\Big(1-\dfrac1k\Big)$ and \\$L^{k^*(p-1),r(p-1)}=\Big(L^{n'(p-1),\infty},L^p\Big)_{\theta,r(p-1)}$.
From the abstract result  Theorem \ref{t2},  we have for $f_1,\ f_2$ in $L^{k,r}(\O)$:
$$||{\mathcal T}f_1-{\mathcal T}f_2||_{\theta,r(p-1);0}\lesssim ||f_1-f_2||^{\frac1{p-1}}_{\theta,r}.$$
Noticing that bounded functions are dense in $L^{k,r}(\O)$, we get the result.
\HF

 Proposition \ref{p3} improves  previous known results considering the case $k=r$: in fact, the usual estimate is only obtained in $\Big[L^{k^*(p-1)}(\O)\Big]^n$ (see for instance \cite{Cianchi-Mazja}).\\

Let us apply now those identifications of the interpolation spaces to obtain precise regularity of the gradient of an entropic-renormalized solution.
\begin{theo}\label{t6}\ \\
Assume (H1) and (H2) and let $m=\dfrac{np}{(n+1)p-n}$ if $2\LEQ p<n$ and $m\in [1,+\infty[$
 if $p\GEQ n,$ $1\LEQ p_2<+\infty,$ $\lambda\in\R$, $0<\theta<1$. Then, the mapping ${\mathcal T}$ is $\dfrac 1{p-1}$-H\"olderian from
$L^{\frac{m'}{m'-\theta},p_2}\Big(\Log L\Big)^\lambda$  into  $L^{p_\theta,p_2(p-1)}\Big(\Log L\Big)^{\frac\lambda{p-1}}$ with $$\dfrac1{p_\theta}=\dfrac{(1-\theta)(n-1)}{n(p-1)}+\dfrac\theta p\qquad\text{and}\qquad m'=\dfrac m{m-1}.$$
If $2\LEQ p<n$ and $\theta=0$,   the mapping ${\mathcal T}$ is $\displaystyle\frac{1}{p-1}-$H\"olderian from $G\Gamma(1,p;t^{-1}(1-Log t)^{\lambda p})$ into $(L^{n'(p-1),\infty}(\Omega), L^p(\Omega))_{0,p(p-1);\lambda/(p-1)}$. This latter space has norm equivalent to
$$
\|f\|_{(L^{n'(p-1),\infty}(\Omega), L^p(\Omega))_{0,p(p-1);\lambda/(p-1)}}\approx \left[ \int_0^1 \left( \sup_{0<s<t^\sigma} s^{\frac{p'}{pn'}}f_*(s)(1-Log t)^{\frac{\lambda}{p-1}}\right)^{p(p-1)}\frac{dt}{t}\right]^{\frac{1}{p(p-1)}},\leqno(*)
$$
where $1/\sigma=p'/(pn')-1/p$.\\
 In the case $\theta=1$, the mapping ${\mathcal T}$ is $\displaystyle\frac{1}{p-1}-$H\"olderian from $(L^{1}(\Omega), L^{(p^*)'}(\Omega))_{1,p;\lambda}$ into\\ $(L^{n'(p-1),\infty}(\Omega), L^p(\Omega))_{1,p/\alpha;\lambda\alpha}$.
The first space has (quasi)norm
	$$\left[ \int_0^1\left(\left(\int_t^1 f_\ast(s)^{(p^*)'}ds\right)^{\frac{1}{(p^*)'}}
			 (1-\log t)^{\lambda }\right)^p \dfrac{dt}{t}\right]^{\frac 1p}$$
			 while the second has (quasi)norm
			 $$\left[ \int_0^1\left(\left(\int_t^1 f_\ast(s)^{p}ds\right)^{\frac{1}{p}}
			 (1-\log t)^{\lambda\alpha}\right)^{\frac{p}{\alpha}} \dfrac{dt}{t}\right]^{\frac{\alpha}{p}}.$$
			 Here $\alpha = \frac 1{p-1}.$
\end{theo}
\Pr
Since ${\mathcal T}$ is $\dfrac 1{p-1}$-H\"olderian from $L^1$  into $[L^{n'(p-1),\infty}(\O)]^n$ and $L^m$ into $[L^p(\O)]^n$, and since smooth functions are dense in Lorentz-Zygmund space $L^{p_1,p_2}(\Log L)^\gamma$, $1\LEQ p_1$, \  $p_2<+\infty$ and $(L^1,L^m)_\tpdl=L^{\frac{m'}{m'-\theta},p_2}(Log L)^\lambda$ according to Proposition \ref{p2.2}, we deduce from Theorem \ref{t2} that ${\mathcal T}$ maps $(L^1,L^m)_\tpdl$ into $\Big(L^{n'(p-1)},L^p\Big)_{\theta,\frac{p_2}\a;\lambda\alpha}$ with $\a=\dfrac1{p-1}$.\\
The identification of the last space   given by Corollary \ref{c1p6} of Proposition \ref{p6}  proves the results.\\
Let $2\LEQ p<n$ and $\theta=0$. First of all, as before, the mapping ${\mathcal T}$ is $\displaystyle\frac{1}{p-1}-$H\"olderian from $L^1(\Omega)$ into $[L^{n'(p-1),\infty}(\Omega)]^n$. On the other hand, in this case, by Lemma 4.3 we know that the mapping ${\mathcal T}$ is $\displaystyle\frac{1}{p-1}-$H\"olderian from $L^{(p^*)'}(\Omega)$ into $[L^p(\Omega)]^n$, where as usual $(p^*)'=np/(np-n+p)$. Noticing that by Proposition \ref{p2.2}, $(L^1,L^{(p^*)'}(\Omega))_{0,p;\lambda}=G\Gamma(1,p;t^{-1}(1-Log t)^{\lambda p})$ and that $L^{(p^*)'}(\Omega)$ is dense therein, we can therefore apply Theorem \ref{t2} and get that ${\mathcal T}$ is $\displaystyle\frac{1}{p-1}-$H\"olderian from $(L^{1}(\Omega), L^{(p^*)'}(\Omega))_{0,p;\lambda}$ into $(L^{n'(p-1),\infty}(\Omega), L^p(\Omega))_{0,p/\alpha;\lambda\alpha}$. Then the assertion follows, because the domain space and the target space have been identified in Proposition 2.1. The same argument holds for $\theta = 1$ .
\HF

To obtain boundedness of the solution in a more general situation, we need to assume  (H3).   We have the following:
\begin{theo}\label{t7}\ \\
Assume  (H1),  (H2) and  (H3). Let $0\LEQ\theta<1,  \ 1<p_2<+\infty,\ \lambda\in\R, \ f\in L^{\frac{n'}{n'-\theta}\,,\,p_2}(\Log L)^\lambda$ $n'=\dfrac n{n-1},\ 0<\theta<1$ and $f\in G\Gamma(1,p_2;w_2)$ with $w_2(t)=t^{-1}\Big(1-\Log t\Big)^{\lambda p_2}$ if $\theta=0$.\\
Then the entropic-renormalized solution  $u$ of the Dirichlet equation (\ref{eq15}) has its gradient in $L^{\frac{n(p-1)}{(1-\theta)(n-1)}\,,\,p_2(p-1)}\Big(\Log L\Big)^{\frac\lambda{p-1}}$ .
\end{theo}
\Pr
Since ${\mathcal T}$ is $\dfrac 1{p-1}$-H\"olderian from $L^1$ into $[L^{n'(p-1,\infty}(\O)]^n$ and ${\mathcal T}$, by Proposition \ref{p1}, is bounded from $L^{n,1}(\O)$ into $L^\infty(\O)$,  then, following Theorem \ref{t1}, ${\mathcal T}$ is bounded from $(L^1,L^n)_\tpdl$
into $\left(L^{n'(p-1)},L^\infty\right)_{\theta,p_2(p-1);\frac\lambda{p-1}}$. With the identification of those interpolation spaces we obtain the result.\HF
\subsection{\bf  Few results  on the case   $1<p<2$}\ \\
Some of the above results remain true in the case $1<p<2$. The fundamental changes concern the H\"older properties than can exist but are not sharp as for the case $p\GEQ2$, and the H\"older constant appearing    depend on the data. Here   is an example.
\begin{theo}\label{t30}{\bf  (local Lipschitz contraction when $1<p<2$)}\ \\
Let $1<p<2$, $p^*=\dfrac{np}{n-p},\ n\GEQ2,$ \
 $\pep=\dfrac{np}{np+p-n}$ its conjugate and $f_1$ (resp $f_2$) in $L^\pep(\O)$. Then, for the weak solution $u$ (resp. $v$) of (\ref{eq12}), say $-\Delta_p u+V(x;u)=f_1$, whenever $V$ satisfies $(H1)$ and $(H2)$, one has:
  \begin{enumerate}
\item $||\Nu||_{L^p}\LEQ c||f_1||^{\frac1{p-1}}_{L^\pep}\,\,$, $\,||\Nv||_{L^p}\LEQ c||f_2||^{\frac1{p-1}}_{L^\pep}\,$.
\item $||\Numv||_{L^p}\LEQ c\Big(||\Nu||_{L^p}+||\Nv||_{L^p}\Big)^{2-p}||f_1-f_2||_{L^\pep}$.
\end{enumerate}
Here the constant $c$ depends only on $p$ and $\O$.
\end{theo}
\Pr
Since we have stability result, we may assume that $f_1$ and $f_2$ are bounded. Arguing as before, one has, using Poincar\'e Sobolev inequality, that
$$\int_\O|\Nu|^pdx+\int_\O uV(x;u)dx\LEQ||f_1||_{L^\pep}\LEQ c||\Nu||_{L^p}\cdot||f_1||_{L^\pep}.$$
Dropping the non negative term  $\DST\int_\O u V(x,u)dx\GEQ0$, we obtain (1).\\
As to the second statement, we use the following inequality (see \cite{LB}  or \cite{DiazNPDE}) concerning the $p$-Laplacian, namely, setting $\WH a (\Nu)=|\Nu|^{p-2}\Nu$, we have
\begin{equation}\label{eq330}\ \\
\Big(\WH a(\Nu)-\WH a(\Nv)\Big)\cdot\Numv\GEQ\a\dfrac{|\Numv|^2}{\Big(|\Nu|+|\Nv|\Big)^{2-p}}\ a.e.\ in\ \O.
\end{equation}
Therefore, making the differences between the two equations and dropping  non negative terms containing $V$, we have from relation (\ref{eq330})
\begin{equation}\label{eq331}
\int_\O\dfrac{|\Numv|^2}{(|\Nu|+|\Nv|)^{2-p}}dx\LEQ c||f-g||_{L^\pep}||\Numv||_{L^p}.
\end{equation}
We have used the Poincar\'e-Sobolev inequality.\\
Now we estimate $\DST\int_\O|\Numv|^pdx.$
Adding the term $\Big(|\Nu|+|\Nv|\Big)^{(p-2)\frac p2}$, the H\"older inequality yields
$$\int_\O|\Numv|^pdx\LEQ\left(\int_\O|\Numv|^2\big(|\Nu|+|\Nv|\Big)^{p-2}dx\right)^{\frac p2}\left[\int_\O\Big(|\Nu|+\Nv|\Big)^pdx\right]^{1-\frac p 2},$$
and with the help of relation (\ref{eq331}), we   have
$$||\Numv||_{L^p}\lesssim ||f_1-f_2||_{L^\pep}\Big(||\Nu||_{L^p}+||\Nv||_{L^p}\Big)^{2-p}.$$
\HF
\begin{coro}[\bf of Theorem \ref{t30}]\label{S4C1}\ \\
Under the same assumptions as in Theorem \ref{t30}, there exists a constant $c$ depending only on $p$ and $\O$ such that
$$||\Numv||_{L^p}\LEQ c\Big(||f_1||^{\frac1{p-1}}_{L^\pep}+||f_2||^{\frac1{p-1}}_{L^\pep}\Big)||f_1-f_2||_{L^\pep}.$$
In particular, the mapping ${\mathcal T}$ is locally Lipschitz  from $L^\pep(\O)$ into $[L^p(\O)]^n$ with $\dfrac1{\pep}=1-\dfrac1p+\dfrac1n$.
\end{coro}
We can have, therefore, the following weaker version of Proposition \ref{p3} when $1<p<2$.

\begin{propo}\label{p31}\ \\
Assume (H1), (H2) and (H3). If $1<p<2$,\ \  $\pep<k<n$, $r\in[1,+\infty]$, 
then the non-linear mapping ${\mathcal T}$ is bounded from $L^{k,r}(\O)$ into $L^{k_1,r}(\O)$, $k_1=\dfrac 1{1-\theta(p-1)}$, with $\theta=p^*(1- \dfrac1 k).$
\end{propo}
\Pr
Following Corollary \ref{S4C1} of Theorem \ref{t30} and Proposition \ref{p1} , the hypotheses of Theorem \ref{t1} are valid for ${\mathcal T}$ with $X_0=L^\pep(\O)$, $X_1=L^{n,1}(\O)$, $\lambda=0$, $Y_0=[L^p(\O)]^n, $ $Y_1= [L^\infty(\O)]^n$, $\a=1$ and $\b=\dfrac 1{p-1}$. According to Theorem \ref{t1}, ${\mathcal T}$ is then a locally bounded mapping from $(X_0,X_1)_{\theta,r}$ into
$(Y_0,Y_1)_{\theta(p-1),r}$ with $\theta\in[0,1]$ such that $\dfrac1k=\dfrac{1-\theta}\pep+\dfrac\theta n$.
Therefore $(X_0,X_1)_{\theta,r}=L^{k,r}(\O)$ and
 $\big(Y_0,Y_1\big)_{\theta(p-1),r}=L^{k_1,r}(\O)\hbox{ with } k_1=\dfrac 1{1-\theta(p-1)}.$
\HF
\begin{rem}\ \\
\begin{enumerate}
\item [a.)] One can make precise the bound for ${\mathcal T}$ locally, according to Theorem \ref{t1}, Corollary \ref{S4C1} of Theorem~\ref{t30}, and Proposition \ref{p1}.
\item [b.)] If $p>n$, then the mapping ${\mathcal T}$ is Lipschitz from $L^1(\O)$ into $[L^p(\O)]^n$: this is a consequence of the Poincar\'e-Sobolev inequality that   we have recalled in Proposition \ref{p1}. Therefore, in view of the Cianchi-Maz'ja's regularity result, the application ${\mathcal T}$ is bounded from $(L^1,L^{n,1})_{\theta,q;\lambda}$ into $[(L^p,L^\infty)_{\theta,q;\lambda}]^n$.
\item [c.)]The list of applications of the above applications is not exhaustive, the reader might derive more results combining those abstract theorems and propositions.
\item [d.)]For other results concerning equations with data in Lorentz spaces, see e.g. \cite{Murat1, Murat2}.
\end{enumerate}
\end{rem}
\section{\bf\large Application of the interpolation for the regularity of the solution of the anisotropic equation}
\subsection{\bf  Preliminary results on anisotropic equations}\ \\
We want to provide similar results as before  for the solution of
\begin{equation}\label{eq100}
\begin{cases}-\Delta_\vecp u+V(x;u)=f&in\ \O\\
u=0&on\ \p\O.
\end{cases}
\end{equation}
Here
$\DST\Delta_\vecp u=-\sum_{i=1}^n\dfrac \p{\p x_i}\left(\left|\dfrac{\p  u}{\p x_i}\right|^{p_i-2}\dfrac{\p u}{\p x_i}\right)$,
$\vecp=(p_1,\ldots,p_n)$, $1<p_i<+\infty$, $\vec {p'}=(p^{'}_1,\ldots,p^{'}_n)$, \\where $p^{'}_i$ is the conjugate of $p_i$.\\
The main differences reside in the exponent appearing in different directions of the
space $\R^n$. Moreover, the estimates concern directly the derivatives in each direction of the $\R^n$-space.\\
Let us recall, from the Introduction, that the real number $p$ is defined as
$\dfrac1p=\dfrac1n\DST\sum_{i=1}^n\dfrac1{p_i}.$ When $\DST\sum_{i=1}^n\dfrac1{p_i}>1$ (say $p<n$) , we set $p^*=\dfrac{np}{n-p} $. We will focus first on the case
$p<n$ for having the H\"olderian property of the mapping ${\mathcal T}$. We set
$$W^{1,\vecp}_0(\O)=\Big\{\f\in W^{1,1}_0(\O)\hbox{ such that }\p_i\f\in L^{p_i}(\O)\ ,\ i=1,\ldots n\Big\}$$
$$\S^{1,\vecp}_0=\Big\{v:\O\to\R\hbox{  measurable  s.t. } \tan^{-1}(v)\in W^{1,1}_0(\O)$$
$$ and \ T_k(v)\in W^{1,\vecp}_0(\O)\hbox{ with }\sup_{k>0}\Big[\Max_{1\LEQ i\LEQ n} k^{\frac1{p_i}}||\p_iT_k(v)||_{L^{p_i}}<+\infty\Big]\Big\}.$$
The definition of an entropic-renormalized solution is similar to Definition \ref{d1111}; we replace the operator and the spaces by the above ones.
\begin{defi}\label{d3100}{\bf entropic renormalized solution for anisotropic equation}\\
For $\xi=(\xi_1,\ldots,\xi_n)\in\R^n$ we consider the vector field
$\WH a _\vecp(\xi)=\Big(|\xi_1|^{p_1-2}\xi_1,\ldots,|\xi_n|^{p_n-2}\xi_n\Big).$\\
We will say that a function $u$ defined on $\O$ is an entropic-renormalized solution associated to the Dirichlet problem
\begin{equation}\label{eq100a}
-\Delta_\vecp u+V(x;u)=f\in\ L^1(\O)\qquad u=0\ on\ \p\O
\end{equation}
if
\begin{enumerate}
\item $u\in\S^{1,\vecp}_0(\O),\ V(\cdot,u)\in L^1(\O).$
\item $\forall\,\eta\in W^{1,\infty}(\O),\ , \forall\,\f\in W^{1,\vecp}_0(\O)\cap L^\infty(\O)$ and all $B\in W^{1,\infty}(\R)$ with $B(0)=0$
\ $B'(\s)=0$ for all $\s$ such that $|\s|\GEQ\s_0>0$, one has:
\end{enumerate}
\begin{equation}\label{eq1150}
\int_\O\WH a_{\vec p}(\Nu).\nabla\Big(\eta B(u-\f)\Big)dx+\int_\O V(x;u)\eta B(u-\f)dx=\int_\O f\eta B(u-\f)dx.
\end{equation}
\end{defi}
Concerning the existence and uniqueness,  let $p^*$ be the number defined for the validity of the Poincar\'e-Sobolev inequality: $\exists c>0$ such that
$$\forall v\in W_0^{1,\vecp}(\O)\qquad\left(\int_\O|v|^{p^*}(x)dx\right)^{\frac1{p^*}}\LEQ c\left(\sum_{i=1}^n\int_\O|\p_i v|^{p_i} dx\right)^{\dfrac1p}.$$
Considering the main operator $$ Au= -\div\big(\WH a_{\vecp}(\Nu)\big)+V(\cdot,u)$$ which is strongly monotonic from $W_0^{1,\vecp}(\O)$ into its dual $W^{-1,\vec {p'}}(\O)$, for $f\in L^1(\O) \cap W^{-1,\vec {p'}}(\O)$, the usual well-known Leray-Lions method or the Leray-Schauder fixed point can be used for having the existence and uniqueness. Moreover, if $f\in L^\infty(\O), $ the maximum principle holds true, using
 for instance the rearrangement technique (see for instance \cite{ABF, DiazNPDE, FPR, Rakobook}) and noticing that the operator $\WH a_{ \vec p}$ satisfies the following coercivity condition:
 there exists  $c_1>0$ such that for all $\xi\in \R^n$, \\
 \C{$\WH a_{\vec p}(\xi)\cdot \xi\GEQ |\xi|^{p_-}-c_1$ with $p_-=\Min(p_i,\ i\in\{1,\ldots ,n\})$.}
  Once the $L^\infty$-estimate is obtained, one may apply standard techniques (approximation method and compactness results) (see
 \cite{Lionslivre, DiazNPDE, Rakotoson1}) to obtain  the following proposition:
\begin{propo}\label{p4000}\ \\
Let  $f\in L^1(\O) \cap W^{-1,\vec {p'}}(\O)$. Then we have a unique weak solution $u\in W^{1,\vecp}_0(\O)$ satisfying:
\begin{equation}\label{eq400}
\int_\O\WH a_\vecp(\Nu)\cdot\nabla\,\f\, dx+\int_\O\f\,V(x;u)dx=\int_\O f\,\f dx,\quad\forall \f\in W^{1,\vecp}_0(\O)\cap L^\infty(\O).
\end{equation}
Moreover, one has the following energy estimates, for $ f \in L^{(p^*)'}(\O),\ p<n, $
\begin{equation}\label{eq401}
\sum_{i=1}^n\int_\O|\p_i u|^{p_i}\,dx+\int_\O u\,V(x;u)\,dx\LEQ c||f||^{p'}_{L^{\pep}},
\end{equation}
where the constant $c$ depends only  on $\O$ and $p$.\\
If $f\in L^\infty(\O)$, then $u\in L^\infty(\O)$, and there are constants $c_i$ independent of $V$ and $f$ so that:
\begin{equation}\label{eq402}
||u||_\infty\LEQ c_1+c_2||f||_\infty^{\frac{p'_-}{p_-} }
\end{equation}
with $p_-=\Min\big(p_i,\ i=1,\cdot ,n\big)$ and $p'_-$ its conjugate.
\end{propo}
\begin{rem}\ \\
A large literature is devoted to the existence for anisotropic equations, besides the above references, one also has  \cite{ABF, FGK, ACCZ}. Those works do not treat the question of local Holderian properties of the gradient as we did here.
\begin{enumerate}
\item [a.)] The fact that the constants $c_1$ and $c_2$ in relation (\ref{eq402}) do not depend on $V$ is due to the hypothesis on $V$ which implies that $\sigma V(x;\s)\GEQ0,\ \ \forall\,\s\in\R$.
\item[b.)] Compactness results concerning anisotropic equation in general form can be found in \cite{ER} (see also \cite{Rako2, Rako4}).
\item[c.)] When $f\in L^\pep(\O)$, the weak formulation is equivalent to the entropic-renormalized formulation. The proof is the same as in \cite{Rako3}.
\item[d.)] The entropic-renormalized solution is specially made for $f\in L^1(\O)$. But the proof of the uniqueness for the solution of (\ref{eq1150})  (see Definition \ref{d3100}) is the same as Benilan et al. \cite {Benilanetal} or Rakotoson \cite{Rako3},   since the operator
$$Au=-{{\rm div\,}}\big(\WH a_\vecp(\Nu)\big)+V(\cdot;u)$$
is monotonic.
 It can be shown that, if $u_1$ and $u_2$ are two solutions in a $T$-space $\S_0^{1,\vecp}(\O)$, then necessarily, one has for all $k>0$
$$\int_{|u_1-u_2|\LEQ k}\Big[\WH a_\vecp (\Nu_1)-\WH a_\vecp(\Nu_2)\Big]\cdot\nabla( u_1-u_2)dx\LEQ0.$$
As to  the existence, it follows using standard approximation technique by replacing $f\in L^1(\O)$ by the sequence $\DST f_j \in L^\infty(\O)$ such that $||f-f_j||_{L^1}\xrightarrow[j\to\infty]{}0, $ $||f_j||_1\LEQ||f||_1$.  Then, one can obtain uniform estimates for the unique weak solution $u_j\in W^{1,\vecp}_0(\O)\cap L^\infty(\O)$
\begin{equation}\label{eq403}
-\Delta_\vecp u_j+V(x;u_j)=f_j.
\end{equation}
\end{enumerate}
\end{rem}
The proof of the following theorem follows the same arguments as in  \cite{Benilanetal} and \cite{Rako2, Rako3, Rako4, Rako5, Rako7}.
\begin{theo}\label{t3101}\ \\
Assume that (H1) and(H2). Then there is a unique entropic-renormalized solution $u$ of (\ref{eq1150}) given in Definition \ref{d3100}. Moreover, for a subsequence denoted by $u_j$, $Du_j(x)\to Du(x)$ a.e. in $\O$.
\end{theo}
\begin{rem}\ \\
In the next paragraph, we will give new and precise spaces where the gradient should be,  under various conditions.   In the case
$$\DST\Min_i p_i=p_{-} \GEQ \Max\Big(\dfrac{p'}{n'};1\Big),\ n'=\dfrac n{n-1},p'=\dfrac p{p-1}\hbox{ conjugate of }p,$$  we have $u\in W_0^{1,1}(\O).$
\end{rem}
\subsection{\bf  The definition of the mappings $\tilde {{\mathcal T}_i} $ from $L^1(\O)$ into $ L^{\frac{n'p_i}{p'},\infty}(\O)$}
\begin{theo}\label{tII21}\ \\
Let $u$ be the entropic-renormalized solution of equation (\ref{eq100}). Then, there exists a constant $c>0$ independent of
$u$ and $f$ such that :
\begin{enumerate}
\item ${\rm meas\,}\{|u|>k\}\LEQ c||f||_{L^1(\O)}^{\frac{p^*}p}k^{-\frac{p^*}{p'}}$,\ \ $\forall\,k>0$.\\
\item $\DST\left\|\dfrac{\p u}{\p x_i}\right\|_{L^{\frac{n'p_i,}{p'}\infty}(\O)}\LEQ c||f||_{L^1(\O)}^{\frac{p'}{p_i}},\qquad\qquad i=1,\ldots,n.$
\end{enumerate}

\end{theo}
\Pr
For the statement  (1), we follow the arguments of Benilan et al \cite{Benilanetal}
  so we drop it. \\
A similar result as for the second statement (2) is given in  \cite{ACCZ}, but the estimate is not precise as we announce here. More, our method is completely different.
To prove it, we apply the  fundamental lemma of Benilan type
(see Lemma \ref{lBt}, in the third paragraph) choosing $h=\left|\dfrac{\p u}{\p x_i}\right|^{p_i}$ and $g=|u|$, to deduce that for
$\lambda>0$ and for all $k>0$:
\begin{equation}\label{eq410}
\meas\left\{\left|\dfrac{\p u}{\p x_i}\right|^{p_i}>\l  \right\}
\LEQ \dfrac 1\l\int_{|u|\LEQ k}\left|\dfrac{\p u}{\p u_i}\right|^{p_i}dx+\meas\{|u|>k\}
\end{equation}
\begin{equation}\label{eq411}
\LEQ\dfrac k\lambda||f||_{L^1}+c_s||f||_{L^1}^{\frac{p^*}p} k^{-{\frac{p^*}{p'}}}.
\end{equation}
This implies
$$\meas\left\{\left|\dfrac{\p u}{\p x_i}\right|^{p_i}>\l\right\}\lesssim
\Min_{k>0}\left(\dfrac1\l||f||_{L^1}k+||f||_{L^1}^{\frac{p^*}p} k^{-{\frac{p^*}{p'}}}\right).$$
Computing the infimum, one has
$$\meas\left\{\left|\dfrac{\p u}{\p x_i}\right|^{p_i}>\l\right\}\lesssim
||f||_{L^1}^{a+1}\l^{-\nppp}\hbox{ with }a=n'\left(\dfrac1p-\dfrac1{p^*}\right).$$
This last inequality implies the result.\HF
In order to derive a H\"olderian mapping, we will use, as in \cite{DiazNPDE,LB}, the elementary inequalities (\ref{eq1300}) and (\ref{eq1301}).

We will deal with different situations. Let us start with the case $p_i\GEQ 2$ for all $i.$
\begin{theo}\label{tII.2.2}\ \\
Assume  (H1), (H2) and  $p_{-} \GEQ \Max\Big(\dfrac{p'}{n'};2\Big)$. Let $i\in\{1,\ldots,n\}$ . Then, the mapping
$$ \WT {\mathcal T}_i : \begin{matrix} L^1(\O)&\longrightarrow&L^{\frac{n'p_i}{p'},\infty}(\O)\\
f&\longmapsto&\WT {\mathcal T}_i f=\dfrac{\p u}{\p x_i}\end{matrix}$$
where $u$ is the unique entropic-renormalized solution is:
\begin{enumerate}
\item $\frac {p'}{p_i}$ -H\"olderian if $p'<p_i,$
\item globally Lipschitz if $p'=p_i$,
\item locally Lipschitz if $p'>p_i $.

\item More precisely, we have a constant $M_1>0$ such that for all $f_1$ and $f_2$ in $L^1(\O)$
$$||\WT {\mathcal T}_i f_1-\WT {\mathcal T}_i f_2||_{L^{\frac{n'p_i}{p'},\infty}}\LEQ M_1||f_1-f_2||_{L^1}^{\frac{p'}{p_i}},\qquad i\in\{1,\ldots,n\}.$$
\end{enumerate}
\end{theo}
\Pr
Let $f_1$ and $f_2$ be in $L^1(\O)$. Due to the stability property, we may assume that $f_1$ and $f_2$ are in $L^\infty(\O)$. Let $u_1$ (resp. $u_2$) be the weak solution of (\ref{eq100}) associated to $f_1$ (resp. $f_2$). Then, for all $k>0$, using relation (\ref{eq100a}) one has:
\begin{equation}\label{eq415}
\alpha\sum_{m=1}^n\int_\O|\p_m T_k(u_1-u_2)|^{p_m}dx\LEQ||f_1-f_2||_{L^1}
\end{equation}
Arguing as in Theorem \ref{tII21}, one deduces that $\forall k>0$
\begin{equation}\label{eq416}\ \\
\meas\{|u_1-u_2|>k\}\LEQ c_\a||f_1-f_2||_{L^1}^{\frac{p^*}p} k^{-{\frac{p^*}{p'}}}.
\end{equation}
From relation (\ref{eq416}), by the same argument as before, which uses   the  fundamental lemma of Benilan type (see Lemma \ref{lBt}, in the second paragraph)  with appropriate choices of $h$ and $g$, 
we deduce
$$||\p_i(u_1-u_2)||_{L^{\frac{n'p_i}{p'},\infty}}\LEQ c||f_1-f_2||_{L^1}^{\frac{p'}{p_i}}.$$
This gives the result.\HF
We have another H\"olderian mapping when the data is in $L^{(p^*)'}(\O)$
\begin{theo}\label{tII.2.3}\ \\
Assume  (H1) and (H2). Let  $\DST\dfrac1p=\dfrac1n\sum_{i=1}\dfrac1{p_i}$
with $\DST \sum_{i=1}^n\dfrac1{p_i}>1$, and let $f_1$ and $f_2$ be two functions $\Lpep(\O)$ with $p^*=\dfrac{np}{n-p} $. Furthermore, we assume  that $p_-\GEQ 2 $. Then, for two weak solutions $u_1$ and $u_2$ associated to $f_1$ and $f_2$, one has
\begin{enumerate}
\item $\DST\sum_{i=1}^n\int_\O|\p_i(u_1-u_2)|^{p_i}dx\LEQ c||f_1-f_2||^{p'}_\Lpep,$
\item $||\p_i(u_1-u_2)||_{L^{p_i}}\LEQ c||f_1-f_2||^{\frac{p'}{p_i}}_{\Lpep}$ for $i=1,\ldots,n$.
\end{enumerate}
\end{theo}
\Pr
The proof is straightforward using $u_1-u_2$ as a test function in the weak formulation for equation (\ref{eq100a}).
${\ }$\HF
Now, we apply the abstract  results concerning interpolations, at first for usual Lorentz spaces as we did before.
\begin{theo}\label{tII.2.4}\ \\
Assume (H1), (H2), and $p_-\GEQ\Max(2;p')$, $1\LEQ k\LEQ (p^*)'$, $r\in [1,+\infty]$. Then for each $i\in\{1,\ldots ,n\}$, the application $\WT {\mathcal T}_i$ is an H\"olderian mapping from  $L^{k,r}(\O)$ into $L^{\frac{k^*p_i}{p'},\frac{rp_i}{p'}}(\O)$, with $k^*=\dfrac{kn}{n-k}$. More precisely, for all $f_1$, $f_2$ in $L^{k,r}(\O)$, $\WT {\mathcal T}_i f_j=\dfrac{\p u_j}{\p x_i}$, $i=1,\ldots,n$, $j=1,2$, we have
$$||\WT {\mathcal T}_if_1-\WT {\mathcal T}_if_2||_{L^{\frac{k^*p_i}{p'},\frac{rp_i}{p'}}}\LEQ M_2||f_1-f_2||^{\frac{p'}{p_i}}_{L^{k,r}}. $$
\end{theo}
\Pr
We argue as in Proposition \ref{p3}, following Theorem \ref{t2}. We have
$$||\WT {\mathcal T}_if_1-\WT {\mathcal T}_if_2||_{(L^{\frac{n'p_i}{p'},\infty} ,L^{p_i})_{\theta,\frac{rp_i}{p'}}}\lesssim ||f_1-f_2||^{\frac{p'}{p_i}}_{(L^1,L^{(p^*)'})_{\theta,r}}$$
whenever $\theta=p^*\left(1-\dfrac1k\right)$, and the identification process (Proposition \ref{p2}) shows that
$$\Big(L^{\frac{n'p_i}{p'},\infty}, L^{p_i}\Big)_{\theta,\frac{rp_i}{p'}}=
L^{\frac{k^*p_i}{p'},\frac{rp_i}{p'}}.$$
This gives the results.
\HF

We may also use the interpolation with a function $(1-Log t)^\lambda$. Here  is an example.
\begin{theo}\label{t461}\ \\
Assume (H1), (H2), and $p'\LEQ p_i$ for each
$i\in\{1,\ldots,n\},\ m=(p^*)',\ \lambda\in\R$ $1\LEQ q_1<+\infty,$ $$\hbox{$\WT {\mathcal T}_i$ is $\dfrac{p'}{p_i}$-H\"olderian mapping from $L^{\frac{p^*}{p^*-1},q_1}
\Big(\Log L\Big)^\lambda$ into $L^{\frac{p_{\theta i},q_1}{\alpha_i} }\Big(\Log L\Big)^{\lambda \alpha_i}$}$$ with
$\dfrac1{p_{\theta i}}=\dfrac{1-\theta}{r_i}+\dfrac\theta{p_i}$, $r_i=\dfrac{n'p_i}{p'}$, $\alpha_i=\dfrac{p'}{p_i}$.
\end{theo}
\Pr
We apply the abstract result stated in  Theorem \ref{t2}, with
$$X_0=L^1,      \quad X_1=L^{(p^*)'},\quad Y_0=L^{r_i,\infty},\quad Y_1=L^{p_i},\quad r_i=\dfrac{n'p_i}{p'},$$
the H\"older exponent being $\alpha_i=\dfrac{p'}{p_i}.$
Since $\WT {\mathcal T}_i$ is $\dfrac{p'}{p_i}$-H\"olderian mapping from $L^1$ into $L^{r_i,\infty}$ and from $L^{(p^*)'}$ into $L^{p_i}$, we deduce that
$$\WT {\mathcal T}_i:\Big(L^1,L^{(p^*)'}\Big)_{\theta,q_1;\lambda }\longrightarrow \Big(L^{r_i,\infty},L^{p_i}\Big)_{\theta,\frac {q_1}{\alpha_i};\lambda \alpha_i}$$
is $\alpha_i$-H\"olderian mapping and the identification space gives the right result.
\HF
We can have similar results for variable exponent but computations are more complicated and are not optimal. So we restrict ourselves to  some estimates.
\subsection{\bf  The Local H\"olderian mappings for the case $\DST\sum_{i=1}^n\dfrac1{p_i}<1 $}
The purpose of this paragraph is to show the following result,  which deals with the case $\DST\sum_{j=1}^n\dfrac1{p_j}=\dfrac np<1.$
\begin{theo}\label{t4301}\ \\
Assume (H1) and (H2). Let $f\in L^1(\O),\ p>n$. Then the unique solution $u\in W_0^{1,\vec p}(\O)$ of the equation (\ref{eq100}) satisfies:
\begin{description}
\item[i.)] $||u||_\infty\LEQ c||f||_1^{\frac1{p-1}}.$\\
\item[ii.)] $\DST\sum_{i=1}^n|| \p_iu||_{p_i}^{p_i}\LEQ c||f||_1^{p'},\quad \dfrac1{p'}+\dfrac1p=1.$\\
\item[iii.)]
In particular if $u_1$ (resp. $u_2$) is the solution associated to $f_1$ (resp. $f_2$), we have for $i\in\big\{1,\ldots,n\big\}$
\begin{enumerate}
\item If $p_i\GEQ2$, then:
$$||\p_i(u_1-u_2)||_{p_i}^{p_i}\lesssim||f_1-f_2||_1\,||u_1-u_2||_\infty.$$
\item If $p_i<2$, then :
$$\int_\O\dfrac{|\p_i(u_1-u_2)|^2}{(|\p_iu_1|+|\p_iu_2|)^{2-p_i}}dx\LEQ||f_1-f_2||_1\,||u_1-u_2||_\infty.$$
\end{enumerate}
\end{description}
\end{theo}
\Pr
Note that when $p>n$, $L^1(\O)$ is a subspace of the dual of $W_0^{1,\vec p}(\O)$, therefore the existence and uniqueness follows from standard theorem concerning monotone operators (see Lions \cite{Lionslivre}) or using
fixed point theorems. So we have for the solution $u\in W_0^{1,\vec p}(\O)$, noticing that $uV(x,u)\GEQ0$, that
\begin{equation}\label{eq4301}
\sum_{i=1}^n||\p_iu||_{p_i}^{p_i}\LEQ\int_\O fudx\LEQ||f||_1\,||u||_\infty.
\end{equation}

Now we use the convexity of the exponential function. Setting temporarily $\lambda_i=\dfrac p{n p_i}$, one has $\displaystyle \sum_{i=1}^n\lambda_i=1$, and setting $a_i=||\partial_i u||_{p_i}$, one has 
$$\left[\prod_{i=1}^n a_i\right]^{p/n}
=e^{ \left[\sum_{i=1}^n \lambda_i\Log a_i^{p_i}\right]}\LEQ
\sum_{i=1}^n \lambda_ia_i^{p_i}\LEQ \sum_{i=1}^n a_i^{p_i}.$$

Hence
\begin{equation}\label{eq4302}
\left[\prod_{i=1}^n||\p_iu||_{p_i}\right]^{\frac1n}\LEQ\left(||f||_1\,||u||_\infty\right)^{\frac1p}.
\end{equation}
Using the Poincar\'e-Sobolev inequality given in Corollary \ref{c1tN1} of Theorem \ref{tlN1}, we derive
\begin{equation}\label{eq4303}
||u||_\infty\lesssim||f||_1^{\frac1{p-1}}.
\end{equation}
Combining relations (\ref{eq4301}) and (\ref{eq4303}), we get the statement i.).\\
Let $u_1$ (resp. $u_2$) be the solution associated to $f_1$ (resp. $f_2$).\\
Since $\Big(V(x;u_1)-V(x;u_2)\Big)(u_1-u_2)\GEQ0$, equation (\ref{eq100}) implies, using elementary inequalities (see relations (\ref{eq1300}) and (\ref{eq1301})), that
$$\sum_{\{i:p_i\GEQ2\}}||\p_i(u_1-u_2)||_{p_i}^{p_i}+\sum_{\{i:p_i <2\}}\int_\O\dfrac{|\p_i(u_1-u_2)|^2}{(|\p_iu_1|+|\p_iu_2)^{2-p_i}}dx\LEQ||f_1-f_2||_1\,||u_1-u_2||_\infty,$$
from which we derive the result.\HF

\begin{coro}[\bf of Theorem \ref{t4301}]\label{c1t4301}
Let $p>n$, $i\in\big\{1,\ldots,n\big\}$. Then the mapping
$\WT {\mathcal T}_i:\begin{matrix}L^1(\O)&\longrightarrow&L^{p_i}(\O)\\
f&\longmapsto&\dfrac{\p u}{\p x_i},\end{matrix}$
where $u$ is the unique solution of (\ref{eq100}), satisfies \begin{description}
\item[$1^{st}$ case]\ \\ If $p_i\GEQ 2$, then $\WT {\mathcal T}$ is a locally $\dfrac1{p_i}$-H\"olderian mapping and for $f_1\in L^1(\O)$, $ f_2\in L^1(\O)$
$$||\WT {\mathcal T}_if_1-\WT {\mathcal T}_i f_2||_{p_i}\lesssim\Big[||f_1||^{\frac1{p-1}}+||f_2||^{\frac1{p-1}}\Big]^{\frac1{p_i}}||f_1-f_2||_1^{\frac1{p_i}}.$$
\item[$2^{nd}$ case]\ \\  If $1<p_i<2,$ then $\WT {\mathcal T}_i$ is a locally $\dfrac12$-H\"olderian mapping and
$$||\WT {\mathcal T}_if_1-\WT {\mathcal T}_i f_2||_{p_i}\lesssim G_0\Big(||f_1||;||f_2||\Big)||f_1-f_2||_1^{\frac12}$$
with $G_0(t;\s)=\Big(t^{p'}+\s^{p'}\Big)^{\frac1{p_i}-\frac12}\Big(t^{\frac1{p-1}}+\s^{\frac1{p-1}}\Big)^{\frac12}$
 for $(t,\s)\in [0,+\infty[\times[0,+\infty[$.
\end{description}
\end{coro}
\Pr
If $i$ is such that $p_i\GEQ2$, then following Theorem \ref{t4301},
$$||\p_i(u_1-u_2)||_{p_i}\lesssim||f_1-f_2||_1\,\Big[||u_1||_1+||u_2||_\infty\Big]\LEQ\left[||f_1||_1^{\frac1{p-1}}+||f_2||_1^{\frac1{p-1}}\right]\,||f_1-f_2||_1.$$
This gives the first statement.\\
Let $i$ be such that $1<p_i<2$. From H\"older's inequality, using Theorem \ref{t4301} iii.), we have
 \begin{equation}\label{eq4304}
||\p_i(u_1-u_2)||_{p_i}^{p_i}\LEQ\Big[||f_1-f_2||_1\,||u_1-u_2||_\infty\Big]^{\frac{p_i}2}\Big[||\p_iu_1||_{p_i}^{p_i}+||\p_i u_2||^{p_i}\Big]^{1-\frac{p_i}2}.
\end{equation}
Using i.) and ii.) of Theorem \ref{t4301},
$$||\p_i(u_1-u_2)||_{p_i}\lesssim\left[||f_1||_1^{\frac1{p-1}}+||f_2||^{\frac1{p-1}}\right]^{\frac12}
\Big[||f_1||_1^{p'} + ||f_2||_1^{p'}\Big]^{\frac1{p_i}-\frac12}||f_1-f_2||_1^{\frac12}.$$
This gives the results.
\HF
As we  observed, if $p_i\GEQ 2\ \forall\,i$, we may have a global-H\"olderian or Lipschitzian mapping.
\begin{coro}[\bf of Theorem \ref{t4301}]\label{c2t4301}
If $\DST p_-=\Min_i p_i\GEQ2$, then for all $i\in\Big\{1,\ldots,n\Big\}$
$$||\WT {\mathcal T}_if_1-\WT {\mathcal T}_if_2||_{p_i}\lesssim||f_1-f_2||_1^{\frac{p'}{p_i}}\quad\forall\, f_1\ and\ \forall\,f_2\ in\ L^1(\O).$$
\end{coro}
\section{\bf \large Few estimates for the solution of $-\Delta_{p(\cdot)}u+V(x;u)=f\in L^1(\O)$}
\subsection{\bf Existence and uniqueness for $-\Delta_{p(\cdot)}u+V(x;u)=f\in L^\infty(\O)$.}\ \\
For the $p(\cdot)$-Laplacian associated to variable exponent, we shall consider the same framework that we introduced in the first paragraph, in particular $p:\O\to]0,+\infty[$, will be a bounded log-H\"older continuous function
$$1<p_-=\DST\Min_{x\in\O} p(x)\LEQ p_+=\Max_{x\in\O} p(x)<n,\ p^*(x)=\dfrac{np(x)}{n-p(x)},$$ whose conjugate is denoted by  $\qq=(p^*)'(\cdot)$. Moreover, we set
$$p'(x)=\dfrac{p(x)}{p(x)-1},\ x\in\OV\O,\ p^*_-=\Min_{x\in\O} p^*(x),\ p_+^*=\Max_{x\in\O} p^*(x),\hbox{ idem for $p'$ conjugate of $p$}.$$
For convenience, we shall add the following assumption for $V$:
$$(H4) : \exists\,\eps>0,\ f_0\in\R_+,\hbox{ such that }\sign(t)V(x;t)\GEQ|t|^\eps-f_0,\hbox{ for a.e. $x\in\O,$ all $t\in\R$}.$$
Such assumption is true if for instance $V(x;t)=|t|^{p(x)-2}t$ with $\eps=p_--1$. We need (H4) only to ensure boundedness of solution when the right hand side is bounded. We first have:
\begin{propo}\label{p5.1}\ \\
Assume (H1), (H2), and (H4), and let $f$ be in $L^\infty(\O)$. Then  we have a unique element $u\in W^{1,p(\cdot)}_0(\O)\cap L^\infty(\O)$ such that:
\begin{equation}\label{eq5001}
\int_\O|\nabla u |^{p(x)-2}\nabla u\cdot\nabla\f\,dx+\int_\O\f V(x;u)dx=\int_\O f\,\f\,dx\qquad\forall\,\f\in W_0^{1,p(\cdot)}(\O).
\end{equation}
Moreover, we have
\begin{equation}\label{eq5002}
\int_\O|\nabla u|^{p(x)}dx+\int_\O uV(x;u)dx\LEQ C\Big[||f||^{p'_-}_{(p^*(\cdot))'}+||f||^{p'_+}_{(p^*(\cdot))'}\Big]
\end{equation}
\begin{equation}\label{eq5003}
||u||_\infty\LEQ M+1,\\
\hbox {with } \Big(f_0+||f||_\infty\Big)^{\frac1\eps} \dot=M.
\end{equation}
\end{propo}
{\bf Idea of the proof}\\
Let $k=M+1$, and define the operator $A$ from $W=W_0^{1,p(\cdot)}(\O)$ into $W'=W^{-1,p'(\cdot)}(\O)$ by
$$Av=-\Delta_{p(.)}v+V(\cdot;T_k(v)),\ v\in W.$$
Due to the assumption (H1) and (H2) on $W$, one can check that $A$ is hemi-continuous, monotonic and coercive (see Lions's book for the definition \cite{Lionslivre}).\\
Therefore, $\forall\,f\in W'$, we have an element $u\in W:Au=f.$
Since the $p(\cdot)$-Laplacian is strictly monotonic and $L^\infty(\O)\subset W'$, we deduce that $u$ is unique and solves
\begin{equation}\label{eq5004}
\int_\O|\nabla u|^{p(x)-2}\nabla u\cdot \nabla\f\,dx+\int_\O\f\,V\big(x;T_k(u)\big)dx=\int_\O f\,\f\,dx,\quad\hbox{ for all }\f\in W_0^{1,p(\cdot)}(\O).
\end{equation}
Let us show  the $L^\infty$-estimates. For this purpose, we consider $$\f=\Big(|T_k(u)|-M\Big)_+\sign (u)\in W_0^{1,p(\cdot)}(\O).$$
Then, dropping the first term, we have:
\begin{equation}\label{eq5005}
\int_\O\Big(|T_k(u)|-M\Big)_+\sign(u)\cdot V\Big(x;T_k(u)\Big)dx\LEQ ||f||_\infty\int_\O\Big(T_k(u)|-M\Big)_+dx.
\end{equation}
Taking into account the hypothesis (H4), we derive from relation (\ref{eq5005}) that
\begin{equation}\label{eq5006}
\int_\O\Big(|T_k(u)|-M\Big)_+\Big[|T_k(u)|^\eps-\Big(f_0+||f||_\infty\Big)\Big]dx\LEQ0.
\end{equation}
The set $\Big\{|T_k(u)|>M\Big\}$ is equal to $\Big\{ |T_k(u)|^\eps>\Big(f_0+||f||_\infty\Big)\Big\}$.
So we deduce from (\ref{eq5006}) that $\Big\{|T_k(u)|>M\Big\}$ is of measure zero, i.e. $|T_k(u)|\LEQ M$ a.e. in $\O$.
But $k>M$ and this implies that $|u(x)|\LEQ k$ almost everywhere in $\O$. This relation and equation (\ref{eq5004})
 imply that $u$ is a solution of (\ref{eq5001}). The uniqueness follows from the fact that
$$\Big[\WH a_{p(\cdot)}(\xi)-\WH a_{p(\cdot)}(\xi')\Big][\xi-\xi']>0\quad\hbox{ if\ }\xi\neq\xi',\ \WH a_{p(\cdot)}(\xi)=|\xi|^{p(x)-2}\xi.\qquad\qquad\diamondsuit$$

\begin{rem}\ \\
\begin{itemize}
\item We may have $||u||_\infty\LEQ M=(f_0+||f||_\infty)^{\frac1\eps}$ if $f_0>0$ or $||f||_\infty>0$, using the same argument but choosing $k=M$,\ $\f=\Big[|T_k(u)|-M+\eta\Big]_+\sign(u)$ with $\eta$ small enough so that $M>\eta$.\\
\item The energy inequality is obtained by choosing $\f=u$ and applying Poincar\'e-Sobolev inequality to derive $$\int_\O fu\LEQ c||f||_{\pointetoile'}||\nabla u||_\point.$$
Using Proposition \ref{p1112}, we have
$$||\Nu||_\point\LEQ\left(\int_\O|\Nu(x)|^{p(x)}dx\right)^{\frac1{p-}}+\left(\int_\O|\Nu(x)|^{p(x)}dx\right)^{\frac1{p+}}, $$
and therefore
$$\int_\O|\Nu|^{p(x)}dx+\int_\O uV(x;u)dx\LEQ c\left[||f||^{p'_-}_{(p^*)'}+||f||^{p'_+}_{(p^*)'}\right].$$
${\ }$ \HF
\item Related existence and uniqueness results are also given in \cite{ BW}. But they do not consider with a lower term and the estimates that we provide here are sharper and precise. More, the compactness provided below is different of their method and we give results on H\"olderian  properties that are not included in their results.
\end{itemize}
\end{rem}

The Proposition \ref{p5.1} is the basis of the existence results when we change the definition of weak solution in (\ref{eq5001}) by entropic solution or renormalized solution, or simply taking the data $f$ in the space $L^1(\O)\cap W^{-1,p'(.)}(\O)$. Here is  an example of such a result:
\begin{corop}[\bf of Proposition \ref{p5.1}]\label{c1p5.1}
For $f\in L^{(\pointetoile)'}$, there exists a unique weak solution $u$ of (\ref{eq5001}) with the test functions $\f\in W_0^{1,\point}(\O)\cap L^\infty(\O),$ which means that
\begin{equation}\label{eq5012}
\int_\O|\Nu|^{p(x)-2}\Nu\cdot\nabla\f\,dx+\int_\O\f\,V(x;u)dx=\int_\O f\f\,dx.
\end{equation}
\end{corop}
{\bf Sketch of the proof}\\
Let $f_j=T_j(f)\in L^\infty(\O)$. Then $\forall\lambda>0 $
$$\int_\O\left|\dfrac{f_j(x)}\lambda\right|^{(p^*(x))'}dx\LEQ\int_\O\left|\dfrac{f(x)}\lambda\right|^{(p^*(x))'}dx.$$
Therefore, $$||f_j||_{[\pointetoile]'}\LEQ||f||_{[\pointetoile]'}.$$
Following Proposition \ref{p5.1}, we have $u_j\in W_0^{1,\point}(\O)$ such that (\ref{eq5001})  and (\ref{eq5002}) hold. We derive
\begin{equation}\label{eq5007}
\int_\O|\Nu_j|^{p(x)}dx+\int_\O u_j V(x;u_j)dx\LEQ c\left[||f||^{p_-}_{[\pointetoile]'}+
||f||^{p_+}_{[\pointetoile]'}\right].
\end{equation}
Since $W_0^{1,\point}(\O)$ is a reflexive space, we have $u\in W_0^{1,\point}(\O)$ and a subsequence still denoted by $u_j$ such that the sequence$ (u_j)_j$ converges weakly to a function $u$ in $W_0^{1,\point}(\O)$, almost everywhere in $\O$ and strongly (by compactness) in $L^{p_-}(\O)$.\\
 Moreover, the fact that $\DST0\LEQ\int_\O u_jV(x;u_j)dx\LEQ C_f<+\infty$ implies that

\begin{equation}\label{eq5008}
\sup_j\int_\O|V(x;u_j)|dx\LEQ C'_f<+\infty.
\end{equation}
Hence we have $\DST\int_\O|V(x;u)|dx\LEQ C'_f$ using Fatou's lemma. Moreover, choosing \\$\f=\Big(|u_j|-t\Big)_+\sign(u_j)$, $t>0$, we derive from (\ref{eq5001})
\begin{equation}\label{eq5009}
\int_{|u_j|>t}|V(x;u_j)|dx\LEQ\int_{|u_j|>t}|f_j|dx.
\end{equation}
Therefore we get
\begin{equation}\label{eq5010}
\lim_{j\to+\infty}\int_\O|V(x;u_j)-V(x;u)|dx=0.
\end{equation}
${\ }$\HF
For the strong convergence of the gradient, we recall the following lemma, which is based on  the monotonicity of the mapping $\WH a_{p(.)} (\xi)=|\xi|^{p(x)-2}\xi$ in our case (see \cite{Rako2, Rako4}).

\begin{lem}\label{l500.1}\ \\
Let $(u_j)_j$ be a sequence of $W_0^{1,\point}(\O)$ having the following properties :
\begin{enumerate}
\item There exists $\qp,\quad 1<q_-\LEQ \qp \LEQ \point$, such that $(u_j)$ remains in a bounded set of $W^{1,\qp}_0(\O)$ and $(u_j)$ converges weakly and a.e. to a function $u$.
\item $z_j^k=T_k(u_j)$ remains in a bounded set of $W_0^{1,\point}(\O)$ for all $k>0.$
\item $\forall\,k>0$, we have a real function $c_k$, such that
$$\forall\,0<\eps<\eps_0,\quad\limsup_{j}\int_{|\,u_j-T_k(u)\,|<\eps}\WH a_\point(\Nu_j)\cdot\nabla\big(u_j-T_k(u)\big)dx\LEQ c_k(\eps)\hbox{ and }\lim_{\eps\to0}c_k(\eps)=0.$$
Then, for a subsequence still denoted by $(u_j)$:
\begin{enumerate}
\item $\Nu_j(x)\DST\xrightarrow[j\to+\infty]{}\Nu(x)\hbox{ a.e in }\O$.
\item If furthermore the conjugate $s$ of $s'(\cdot)\!=\!\dfrac {\qp}{p(\cdot)-1}$ satisfies $\DST\lim_{m \rightarrow \infty}\frac{1}{m}\Big[\int_\O\!s^m(x)dx\Big]^{\frac{1}{m}}\!=\!0$, then
$$\lim_j \int_\O\Big| |\Nu_j|^{q(x)-2}\Nu_j\,dx- \int_\O|\Nu|^{q(x)-2}\Nu\Big|\,dx = 0.$$
\item In particular, for all $\f\in W_0^{1,q'(\cdot)}(\O)$
$$\lim_j\int_\O|\Nu_j|^{q(x)-2}\Nu_j\cdot\nabla\f\,dx=\int_\O|\Nu|^{q(x)-2}\Nu\cdot\nabla\f\,dx.$$
\end{enumerate}
\end{enumerate}
\end{lem}
\Pr
The proof of the first statement is  similar to Lemma 2 of \cite{Rako4} (see also \cite{Rako2}) or Lemma A.5 of \cite{Rako7} for a more general case, so we drop it. But for the second statement, we need to use Theorem 2.1
of  \cite{FGNR} and Vitali's convergence lemma. Indeed, let us set
$$\DST h_j=\Big| |\Nu_j|^{q(x)-2}\Nu_j\,dx- \int_\O|\Nu|^{q(x)-2}\Nu\Big|.$$
Then, the preceding statement shows that $h_j(x) \rightarrow 0$
almost everywhere in $\O$. Besides applying H\"older's inequality, we  have the following uniform integrability, for all measurable set $E \subset \O$:
$$\DST \int_E | h_j(x)|dx \LEQ c||\chi_E||_{s(\cdot)}, $$ for some constant $c$
independent of $j$ and $E$.
Following Theorem 2.1 of  \cite{FGNR}, the condition on $s$ implies that
$||\chi_E||_{s(\cdot)}$ tends to zero as $\meas (E)$ tends to zero. Thus the conditions of Vitali's convergence lemma are fulfilled, so that
$\DST \lim_j \int_\O | h_j(x)|dx =0$ .
\HF
\\
Since $T_k(u) \in W_0^{1,p(\cdot)}(\O)$, for all $k>0$, then for any $\eps >0$
the function $T_\eps(u_j -T_k(u))$ is a suitable test function in relation (\ref{eq5001}).
We then derive the third statement of Lemma \ref{l500.1}. Therefore we have the
necessary convergences for the gradient to pass to the limit in the equation
\begin{equation}\label{eq5011}
\int_\O|\Nu_j|^{p(x)-2}\Nu_j\nabla\f\,dx+\int_\O \f V(x;u_j)dx=\int_\O f_j\f\quad\f\in W_0^{1,\point}(\O)\cap L^{\infty}(\O),
\end{equation}
so that $u$ solves the equation (\ref{eq5012}).\\ The uniqueness is a consequence of strong monotonicity of the $\WH a_\point.$
\HF

\begin{corop}[\bf of Proposition \ref{p5.1},   Local H\"olderian mapping]\label{c2p51} 
Assume that $p_-\GEQ2$. \\The mapping ${\mathcal T}^*:\begin{matrix}L^{(\pointetoile)'}(\O)\longrightarrow \Big[L^\point(\O)\Big]^n\\f\longmapsto {\mathcal T}^*f=\nabla u\end{matrix}$
 is $\a_1=\dfrac{p'_-}{p_+}$-local H\"olderian mapping.\\ Here $u$ is the solution of (\ref{eq5012}) associated to $f$. \\
More precisely, we have: $\forall\,f_1,\ f_2$ in $L^{[\pointetoile]'}(\O)$
$$||{\mathcal T}^*f_1-{\mathcal T}^*f_2||_\point\lesssim\Phi_2(f_1;f_2)||f_1 -f_2||^{\a_1}_\qq$$
where $\Phi(f_1;f_2)=||f_1-f_2||^{1-\a_1}_{L^{[\pointetoile]'}}+||f_1-f_2||^{\a_2-\a_1}_{L^{[\pointetoile]'}}+1$ with $\a_2=\dfrac{p'_+}{p_-}.$
\end{corop}
\Pr
If $u_1$ (resp $u_2$) is the solution of (\ref{eq5001}) with $f=f_1$ (resp $f_2\in L^\qq(\O)$), then$$\int_\O|\Nu_1-\Nu_2|^{p(x)}dx\LEQ c\left[||f_1-f_2||^{p_-}_\qq+||f_1-f_2||^{p_+}_\qq\right].$$
Since we have
$$2||\nabla(u_1-u_2)||_\point\LEQ\left[\int_\O|\nabla(u_1-u_2)^{p(x)}dx\right]^{\frac1{p_-}}+\left[\int_\O|\nabla(u_1-u_2)^{p(x)}dx\right]^{\frac1{p_+}}, $$
 we get the result, noticing that $1\GEQ \a_2=\dfrac{p'_+}{p_-}\GEQ \dfrac{p'_-}{p_+}=\a_1.$
\HF
\subsection{\bf A priori estimates for variable exponents with data in $L^1(\O)$.}\ \\
We only give a priori estimates starting with the equation (\ref{eq5001}).
\begin{propo}\label{p5.2}\ \\
For a solution $u$ of (\ref{eq5001}), one has:
\begin{enumerate}
\item $\DST\int_\O|\nabla T_k(u)|^{p(x)}dx\LEQ k||f||_1,\quad \forall\,k>0.$
\item $\DST||\nabla T_k(u)||_\point\LEQ \Max\Big((k||f||_1)^{\frac1{p_+}};(k||f||_1)^{\frac1{p_-}}\Big)$.
\item $||T_k(u)||_{p^*(\cdot)}\lesssim \Max \Big((k||f||_1)^{\frac1{p_+}};(k||f||_1)^{\frac1{p_-}}\Big)$.
\end{enumerate}
\end{propo}
\Pr
Taking as a test function $\f=T_k(u)$,   we get (1). In order to get (2), we use the estimate
$$2||\nabla T_k(u)||_\point\LEQ\left(\int_\O\Big|\nabla T_k(u)\Big|^{p(x)}dx\right)^{\frac1{p_+}}+\left(\int_\O\Big|\nabla T_k(u)\Big|^{p(x)}dx\right)^{\frac1{p_-}}$$
and statement (1).\\
Finally, the last statement is a consequence of the Poincar\'e-Sobolev inequality.
\HF
Next we want to study the decay of $\meas\Big\{|\Nu|^\point>\lambda\Big\}$, for $\lambda >0$.
To make  our computation easier, we will take the new variable $k=\OV k||f||_1,\ \OV k>0.$
We have:
\begin{propo}\label{p5.4}\ \\
For all $\lambda>0$, all $\OV k>0$, we have
$$\meas\Big\{|\Nu|^\point>\lambda\Big\}\LEQ\dfrac k\lambda+\meas\Big\{|u|>\OV k\Big\}.$$
\end{propo}
\Pr
We use first the fundamental lemma of Benilan type,  see Lemma \ref{lBt}  with $h=|\Nu|^\point$,\  $g=|u|$ and then we apply the statement (1) of the preceding Proposition \ref{p5.2} to conclude.
\HF
Next, we need to estimate the decay of $\meas\Big\{|u|>\OV k\Big\}$. One has:
\begin{propo}\label{p5.5}\ \\
Let $a_1=\dfrac{p^*_+}{p_-}-p^*_-,\ \ \psi_1(t)=\Max\Big(t^{p^*_+};t^{p^*_-}\Big)$
for $t>0$. Assume that $a_1<0$, that is
$$\dfrac{n-p_-}{n-p_+}<p_-\left(\dfrac{p_-}{p_+}\right).$$
Then
$$\meas\Big\{|u|>\OV k\Big\}\lesssim\psi_1\Big(||f||_1\Big)k^{-|a_1|}\quad\hbox{ with } k=\OV k||f||_1.$$
\end{propo}
\Pr
We know that for $\eps<\OV k$, one has $\Big\{|u|>\eps\Big\}=\Big\{|T_{\OV k}(u)|>\eps\Big\}.$ The same argument as before leads to
\begin{equation}\label{eq5200}
\meas\Big\{|u|>\OV k\Big\}\LEQ\Max\left(\dfrac{||f||_1^{p^*_+}}{k^{p^*_+}};
\dfrac{||f||_1^{p^*_-}}{k^{p^*_-}}\right)\int_\O|T_{\OV k}(u)|^{p^*(x)}dx,
\end{equation}
from which we get, using statement (3) of Proposition \ref{p5.2},
$$\meas\Big\{|u|>\OV k\Big\}\lesssim\psi_1\big(||f||_1\big)\Max\Big(k^{-p^*_+};k^{-p^*_-}\Big)\Max\Big(M_1(k)^{p^*_+};M_1(k)^{p^*_-}\Big)$$
where $M_1(k)=\Max\Big(k^{\frac1{p_+}};k^{\frac1{p_-}}\Big).$\\
If $k\GEQ1$, then the above estimate is reduced to
$$\meas\Big\{|u|>\OV k\Big\}\LEQ\psi(||f||_1)k^{-|a_1|},\quad a_1=\dfrac{p^*_+}{p^*_-}-p^*_-.$$
If $k\LEQ1$, then it is reduced to
$$\meas\Big\{|u|>\OV k\Big\}\LEQ\psi_1(||f||_1)k^{a_2},\quad \hbox{with } a_2=\dfrac{p^*_-}{p_+}-p^*_+.$$
But we have
$$a_1-a_2=\dfrac{n^2(p_+-p_-)}{(n-p_+)(n-p_-)}\left[\dfrac1{p_+}+\frac1{p_-}+\frac{n-1}n\right]>0:a_1\GEQ a_2,$$
and therefore for $k\LEQ 1$, $k^{a_2}\LEQ k^{-|a_1|}$.\\
So for all $k>0$, one has \\
\centerline{$\meas\Big\{|u|>\OV k\Big\}\lesssim\psi_1(||f||)k^{-|a_1|}.$}
${\ }$\HF
\begin{theo}\label{t5.1}{(\bf  main estimate for the $L^1$-data)}\\
Under the same assumptions as for Proposition \ref{p5.5},
 there exists a constant $c>0$ depending only on $p,\ n,\ \O$ such that
$$\meas\Big\{|\Nu|^\point>\lambda\Big\}\LEQ c\,\psi_1(||f||_1)^{\frac1{1+|a_1|}}\lambda^{-\frac{|a_1|}{1+|a_1|}}\quad\forall\,\lambda>0.$$
\end{theo}
\Pr
From Proposition \ref{p5.4} and Proposition \ref{p5.5}, we have, for all $k>0$,
$$\meas\Big\{|\Nu|^\point>\lambda\Big\}\LEQ\dfrac k\lambda+c_1\psi_1(||f||_1)k^{-|a_1|}$$
where $c_1$ depends only the Sobolev constant that is on $\O$, $n$, $p$. Taking the infimum of the right hand side, we derive the result.
\HF
\begin{coro}[\bf of Theorem \ref{t5.1}]\label{cor1t5.1} 
Assume that $\dfrac{|a_1|}{1+|a_1|}p_->1$. Then for all $q\in\left[\dfrac{p_+}{p_-},\dfrac{|a_1|}{1+|a_1|}p_+\right[$ we have
$$\int_\O|\Nu|^{\frac q{p_+}p(x)}dx\LEQ c\, \psi_1(||f||_1)^{\frac q{p_+|a_1|}} $$
where $c$ depends only on $\O,\ p,\ n$.
\end{coro}
\Pr
From Theorem \ref{t5.1}, we deduce
$$\int_\O|\Nu|^{\frac q{p_+}p(x)}dx\LEQ c\, \psi_1(||f||)^{\frac q{p_+|a_1|}}\int_0^{|\O|}t^{-\frac{1+|a_1|}{|a_1}\frac q{p_+}}dt<+\infty.$$
\HF

\begin{rem}\ \\
We recover all the condition that we obtained in the preceding section when $p(x)=p$ is constant. In particular, the condition $\dfrac{|a_1|}{1+|a_1|}p_->1$ is equivalent to $p>2-\dfrac1n$ since we have $p\dfrac{|a_1|}{1+|a_1|}=\dfrac n{n-1}(p-1).$
\end{rem}

\subsection{\bf Appendix : An existence and uniqueness result of an entropic-renormalized solution for variable exponents}\ \\
Although it is not the purpose of our paper, we will show now how to prove the existence of an entropic-renormalized solution. The principle is the same as we did in our previous papers, but for convenience, here we give   the main steps.
\begin{theo}\label{t630s}\ \\
Let $q$ be as in Corollary \ref{c1p5.1} of Theorem \ref{t5.1}. Assume (H1), (H2), and (H4), that\\ $q>\dfrac{p_+}{p_-}(p_+-1)$, and let $f\in L^1(\O)$. 
Then there exists a unique solution $u\in W^{1,\qp}_0(\O)$ with $q(x)=\dfrac q{p_+}p(x)$ such that\\
$\forall\,\eta\in W^{1,\infty}(\O),\ \forall\,B\in W^{1,\infty}(\R)\hbox{ with }B(0)=0,\ B'(\s)=0$ for $|\s|\GEQ \s_0>0,$\\ $\forall\,\f\in W_0^{1,\pp}(\O)\cap L^\infty(\O)$
$$\int_\O\WH a_\pp(\Nu)\cdot\nabla\Big(\eta B(u-\f)\Big)dx
+\int_\O\eta B(u-\f)\cdot V(x;u)dx=\int_\O\eta B(u-\f)fdx.$$
\end{theo}
\Pr
We only give the main steps for the existence. Consider $f_j=T_j(f)\in L^\infty(\O)$. Following Proposition \ref{p5.1}, we have  a unique function $u_j\in W_0^{1,\pp}(\O)\cap L^\infty(\O)$ satisfying relation (\ref{eq5001}).
Moreover the above Corollary \ref{c1p5.1} of Theorem \ref{t5.1} shows that $u_j$ remains in a bounded set of $W_0^{1,\qp}(\O)$, and  we have
\begin{equation}\label{eq6301}
\sup_j\int_\O\|\Nu_j|^{q(x)}dx\LEQ c\,\psi_1\big(||f||_1\big)^{\frac q{p_+|a_1|}}.
\end{equation}
Taking as a test function $T_k(u_j)=\f$ in relation (\ref{eq5001}),
we  deduce
\begin{equation}\label{eq6302}
\int_\O|\nabla T_k(u_j)|^{p(x)}dx\LEQ k||f||_1.
\end{equation}
Since $1<\qp<p_+<+\infty$,   the space $W_0^{1,\qp}(\O)$ is reflexive, and we may subtract a sequence still denoted  $u_j$, and have an element $u\in W_0^{1,\qp}(\O)$
 such that
 \begin{itemize}
\item $u_j$ converges weakly to $u$ in $W_0^{1,\qp}(\O)$.
\item$u_j(x)\DST \xrightarrow [j\to +\infty]{}u(x)$ a.e in $\O$.
\item $T_k(u_j)$ converges weakly to $T_k(u)$ in $W_0^{1,p}(\O)$ for all $k>0$.
\end{itemize}
Taking as a test function
$\f=\Big(|u_j|-t\Big)_+\sign(u_j),\  t>0$, and dropping non negative term, we have
\begin{equation}\label{eq6303}
\int_{|u_j|>t}|V(x;u_j)|dx\LEQ\int_{|u_j|>t}|f|dx.
\end{equation}
This relation with the pointwise convergence and assumptions (H1) and (H2), implies
\begin{equation}\label{eq6304}
\lim_j\int_\O|V(x;u_j)-V(x;u)|dx=0.
\end{equation}
Next, we choose as a test function $\f=T_\eps\Big(u_j-T_k(u)\Big)$
with $\eps>0$, so we have
\begin{equation}\label{eq6305}
\int_{|u_j-T_k(u)|<\eps}\WH a_\pp(\Nu_j)\cdot\nabla\big(u_j-T_k(u)\big)dx\LEQ\eps\left[||f||_1+\int_\O|V(x;u_j|dx\right].
\end{equation}
Therefore, we have
\begin{equation}\label{eq6306}
\limsup_j\int_{|u_j-T_k(u)|<\eps}\WH a_\pp(\Nu_j)\cdot\nabla\Big(u_j -T_k(u)\Big)dx\LEQ\eps \left[||f||_1+\int_\O|V(x;u)|dx\right] .
\end{equation}
We may invoke Lemma \ref{l500.1} to derive, for a sequence still denoted $(u_j)$, that $\Nu_j(x)\DST \xrightarrow[j\to +\infty]{}\Nu(x)\hbox{ a.e. in } \O. $
The condition that $q>\dfrac{p_+}{p_-}(p_+-1)$ implies, for all $x$, \\
$q(x)=\dfrac q{p_+}p(x)\GEQ \dfrac q{p_{+}}p_{-}>p_{+}-1\GEQ p(x)-1>p_--1>0.$ Therefore
\begin{equation}\label{eq6307}
\lim_j\int_\O\Big|\,|\Nu_j(x)|^{p(x)-2}\Nu_j(x)-|\Nu(x)|^{p(x)-2}\Nu(x)\Big|dx\equiv0.
\end{equation}
Indeed, let us set $g_j(x)=\Big| \,|\Nu_j|^{\pp-2}\Nu_j-|\Nu|^{\pp-2}\Nu\Big|(x).$ Since $r(x)\dot= \dfrac{q(x)}{p(x)-1}\GEQ\dfrac{qp_-}{p_+(p_+-1)}>1$, \ $r\in C(\OV \O)$, we may apply  \cite[Theorem 2.1]{FGNR} to derive that  for all measurable set $E\subset\O$,
\begin{equation}\label{eq6308}
\lim_{|E|\to0}||\chi_E||_{r'(\cdot)}=0 
\end{equation}
where $r'(x)=\dfrac{r(x)}{r(x)-1}$, and
$\chi_E$ is the characteristic function of  $E$. But the boundedness of the sequence $(u_j)_j$ in $W_0^{1,\qp}(\O)$ and H\"older inequality imply, for all measurable set  $E$, that
\begin{equation}\label{eq6309}
\sup_j\int_E|g_j(x)|dx\LEQ c||\chi_E||_{r'(\cdot)}.
\end{equation}
Thus, we may apply Vitali's convergence theorem to derive
$$\lim_{j\to+\infty}\int_\O|g_j(x)|dx=0$$
since $\DST g_j(x)\xrightarrow[j\to+\infty]{}0$ a.e., so that  we have the uniform integrability given by (\ref{eq6309}).\\
The convergences given by relation (\ref{eq6304}) and relation ({\ref{eq6307}) are enough to prove the existence of a weak solution  when $f\in L^1(\O).$\\
To obtain an entropic-renormalized solution, we need further estimates:
\begin{lem}\label{l6301}{\bf (Gradient behavior)}\ \\
One has for all $m\GEQ 0$, all $ j\GEQ0$:
\begin{enumerate}
\item $\DST\int_{\{x:m\LEQ|u_j|\LEQ m+1\}} |\Nu_j(x)|^{p(x)}dx\LEQ\int_\O|f|\,\Big|T_{m+1}(u_j)-T_m(u_j)\Big|dx.$
\item $\DST\int_{\{x:m\LEQ|u_j|\LEQ m+1\}} |\Nu(x)|^{p(x)}dx\LEQ\limsup_j\DST
\int_{\{x:m\LEQ|u_j|\LEQ m+1\}} |\Nu_j(x)|^{p(x)}dx$$$\DST\LEQ\int_\O|f(x)|\,\Big|T_{m+1}(u)-T_m(u)\Big|dx\xrightarrow[m\to+\infty]{}0.$$
\end{enumerate}
\end{lem}
}
\Pr
We can take as test function $\psi_{mj}=T_{m+1}(u_j)-T_m(u_j)$. Since
$$\int_\O\psi_{mj}\,V(x;u_j)dx\GEQ0,$$
and
$$\int_\O\WH a_\pp(\Nu_j)\cdot\nabla\psi_{mj}dx=\int_{m\LEQ|u_j|\LEQ m+1}|\Nu_j(x)|^\pp dx,$$
we get statement (1). On the other hand,
  statement (2) follows from (1) using Fatou's lemma and pointwise convergences of the gradient for the lower bound and the pointwise convergence of $u_j$ for the upper bound, combined with the Lebesgue dominated convergence.
\HF
For convenience for the next results,  for $v\in L^1(\O)$, we shall denote  $v^m=T_m(v)$ and we define $h_m \in W^{1,\infty}(\R) $ :
$$h_m(\s)=\begin{cases}1&\si\ |\s|\LEQ m,\\
0&\si\ |\s|\GEQ m+1,\\m+1-|\s|&otherwise.\end{cases}$$
\begin{lem}\label{l6302}\ \\
Let $\eta\in W^{1,r}(\O),\ \ r>n,\ \ b\in W^{1,\infty}(\R)\with B(0)=0,\quad B'(\s)=0$ for $|\s|\GEQ \s_0>0,$\\ $\f\in W^{1,\pp}_0(\O)\cap L^\infty(\O)$ and set $\f_{mj}=\psimT \ h_m(u_j)$. Then
\begin{enumerate}
\item $\f_{mj}\in W_0^{1,\pp}(\O)\cap L^\infty(\O),\ \ \forall\, m\GEQ0,\quad\forall\,j\GEQ0$.
\ \\
\item $\DST\left|\int_\O|\Nu_j^{m+1}|^{p(x)-2}\Nu_j^{m+1}\nabla\Big(\psimu\Big)h_m(u_j)+\int_\O\f_{mj}\Big[V(x;u_j)-f_j\Big]dx\right|$$$\LEQ||\eta||_\infty||\,||B||_\infty\int_{m\LEQ |u_j|\LEQ m+1}|\Nu_j|^{p(x)}dx.$$
\end{enumerate}
\end{lem}
\Pr
Taking $\varphi_{mj}$ as a test function, in the relation (64) satisfied by the solution $u_j$, we have :\\
$\displaystyle\left|\int_\Omega |\nabla u_j^{m+1}|^{p(x)-2}\nabla u_j^{m+1}\nabla \Big(\eta B(u^{m+1} - \varphi) \Big)h_m(u_j)+\int_\Omega \varphi_{mj}\Big[V(x;u_j)-f_j\Big]dx\right|$$$= \displaystyle\left|\int_\Omega |\nabla u_j^{m+1}|^{p(x)-2}\nabla u_j^{m+1}\Big(\eta B(u^{m+1} - \varphi) \Big)\nabla h_m(u_j) \right|= A.$$

Since $h_m(u_j)\in W^{1,\pp}(\O)$, and $$|\nabla h_m(u_j)|\LEQ\begin{cases}|\Nu_j|&\si m\LEQ |u_j|\LEQ m+1,\\0& \ewe,\end{cases}$$
the last quantity  $A$ can be estimated as:
$$ A \LEQ||\eta||_\infty||\,||B||_\infty\int_{m\LEQ |u_j|\LEQ m+1}|\nabla u_j|^{p(x)}dx,$$
so we derive statement (2).\\
Let us note that $\psimu$ is in $W^{1,\pp}_0(\O),\ \WH a_p(0)=0$.
\HF
\begin{lem}\label{l6303}\ \\
For fixed $m$,
$\WH a_\pp(\Nu_j^{m+1})\hbox{ converges weakly to $\WH a_\pp(\Nu^{m+1})$ in }[L^{p'(\cdot)}(\O)]^n.$
\end{lem}
\Pr
The pointwise convergence of the gradient implies
$$\WH a_\pp(\Nu_j^{m+1})\xrightarrow[]{}\WH a_p(\Nu^{m+1})\aein \O.$$
Furthermore, we know that
  $$\DST||\WH a_\pp(\Nu_j^{m+1})||_{L^{p'(\cdot)}}\LEQ c_m<+\infty.$$
  By the reflexivity of $[L^{p'(\cdot)}(\O)]^n$, we derive the result.
\HF
\begin{coro}[\bf of Lemma \ref{l6301},  \ref{l6302},  \ref{l6303}]\label{cl630123}\ \\
The function $u$ satisfies, for all $m\GEQ0$
$$\Big| \int_\O\WH a_\pp(\Nu^\mpu)\cdot\nabla(\psimu) h_m(u)+\int_\O h_m(u)\psimu\Big[V(x;u)-f\Big]dx \Big| $$$$\LEQ||\eta||_\infty||B||_\infty\int_\O|f(x)|\,\Big|T_\mpu(u)-T_m(u)\Big|dx.$$
\end{coro}
\Pr
Since $\nabla\Big(\psimu\Big) h_m(u_j)$ converges strongly to $\nabla(\psimu) h_m(u)$ in $[L^\pp(\O)]^n$, combining with  the weak convergence of Lemma \ref{l6303}, we obtain
\begin{equation}\label{eq6310}
\lim_{j\to+\infty}\int_\O\WH a_\pp(\Nu_j^\mpu)\nabla\Big(\psimu \Big)h_m(u_j) dx =\int_\O\WH a_p(\Nu^\mpu)\cdot\nabla(\psimu) h_m(u_j)dx.
\end{equation}
Since $V(\cdot;u_j)$ (resp. $f_j$) converges strongly to $V(\cdot u)$ (resp $f$) in $L^1(\O)$, we have
\begin{equation}\label{eq6311}
\lim_{j\to+\infty}\int\f_{mj}\Big[V(x;u_j)-f_j\Big]dx=\int_\O h_m(u)\psimu\Big[V(x;u)-f\Big]dx.
\end{equation}
Combining with Lemma \ref{l6301}, the two last relations and Lemma \ref{l6302} give the result.
\HF
We then have:
\begin{lem}\label{l6304}\ \\
\vspace{-1cm}
\begin{enumerate}
\item $\DST\lim_{m\to+\infty}\int_\O\WH a_\pp (\Nu^\mpu)\nabla\Big(\psimu\Big) h_m(u)dx=\int_\O\WH a_\pp(\Nu)\cdot \nabla\Big(\eta\,B(u-\f)\Big)dx.$
\ \\
\item $\DST\lim_{m\to+\infty}\int_\O h_m(u)\psimu\Big[V(x;u)-f\Big]dx
=\int_\O\eta B(u-\f)\Big[V(x;u)-f\Big]dx.$
\end{enumerate}
\end{lem}
\Pr
As we have already observed before,
 \begin{equation}\label{eq6311a}
\WH a_\pp(\Nu^\mpu)\cdot\nabla\Big(\psimu\Big)h_m(u)=\WH a_\pp(\Nu)\cdot\nabla\Big(\psiuf\Big)\,h_m(u),
\end{equation}
because of the definition of $h_m$, $\WH a_\pp(0)=0.$\\
Moreover, when expanding the gradient, we have:
$$\WH a_\pp(\Nu)\cdot\nabla\Big(\eta B\big(u -\varphi\big)\Big)=\WH a_\pp(\Nu)\cdot\nabla\eta B(u-\f)+\WH a_\pp(\Nu)\cdot\Nu B'(u-\f)\eta.$$
Since $B'(u-\f)=0$ if $|u-\f|>\s_0$, then, setting $k_0=||f||_\infty+\s_0$,  we  have:
$$\WH a_\pp(\Nu)\cdot\Nu B'(u-\f)\eta=\WH a_\pp(\Nu^{k_0})\cdot\Nu^{k_0}B'(u-\f)\eta.$$
Hence we deduce, from the preceding decomposition,  the following estimate: $$\Big|\WH a_\pp(\Nu)\nabla \Big(\psiuf\Big) \Big| \LEQ c\Big[|\Nu(x)|^{p(x)-1}+|\Nu^{k_0}|^{p(x)}\Big]\dot=R(x).$$
Here  $c>0$ is independent of $u, \f$ .

One has $R\in L^1(\O)$. Therefore, by the Lebesgue dominated theorem, we have:
\begin{equation}\label{eq6312}
 \lim_{m\to+\infty}\int_\O\WH a_\pp(\Nu)\cdot\Big(\nabla \psiuf\Big)h_m(u)=\int_\O\WH a_\pp(\Nu)\nabla\Big(\psiuf\Big).
\end{equation}
Both relations (\ref{eq6311a}) and (\ref{eq6312}) infer the first statement (1) of Lemma \ref{l6304} while the second one comes from the Lebesgue dominated theorem.
\HF
{\bf End of the proof of the main theorem}\\
Letting $m\to+\infty$ in Corollary \ref{cl630123} of Lemmas \ref{l6301} to \ref{l6303} with the help of Lemma \ref{l6304}, we get that $u$ is an entropic-renormalized solution.\\
For the uniqueness, we may use the method of Benilan et al \cite{Benilanetal} since any entropic-renormalized solution is also an entropic solution in their sense. Note that here, in our case, the solution is always in $W^{1,1}_0(\O)$. The second method consists in noticing that since $f_1$ (resp $f_2$) are two elements of $L^1(\O)$ and
$u_1$ (resp $u_2$), we have:
\begin{lem}\label{lem730123}
\begin{equation}\label{eq6313}
\int_\O\dfrac{\Delta[u_1-u_2]}{1+|u_1-u_2|}dx\LEQ\dfrac\pi2\int_\O|f_1-f_2|dx
\end{equation}
\hbox{whenever }$\Delta[u_1;u_2]=\Big[\WH a_\pp(\Nu_1)-\WH a_\pp(\Nu_2)\Big]\cdot\Nu(u_1-u_2)\GEQ 0.$
\end{lem}
The proof of this lemma needs the following result, which can be carried out in an even more general situation:
\begin{lem}\label{lem730124}\ \\
      Let $w$ be in $W^{1,1}_{loc}(\O)$ such that for all $k>0$, $T_k(w)=w^k \in W_0^{1,p(.)}(\O)$ and let $B\in W^{1,\infty}(\R)$ with $B(0)=0$, $B'(\s)=0$ for all $\s$ such $|\s|\GEQ \s_0>0$, $\f \in W_0^{1,p(.)}(\O)\cap L^{\infty}(\O) $. Then $B(w-\f)$ is in $W_0^{1,p(.)}(\O)\cap L^{\infty}(\O) $.
\end{lem}
 {\bf Proof of Lemma \ref{lem730124}}\ \\
If we choose $k=\s_0 + \||\f||_\infty$, then $B(w^k - \f)$ is in
 $W_0^{1,p(.)}(\O)\cap L^{\infty}(\O) $. Moreover, almost everywhere in
 $\O$,
 $$\nabla B(w^k - \f)= B'(w-\f)\nabla (w-\f)= \nabla B(w-\f). $$
 Since $\nabla B(w-\f) \in W^{1,1}_{loc}(\O)$, then  the above equality holds in the sense of distribution and implies the result.
 \HF

 {\bf Proof of Lemma \ref{lem730123}}\ \\
 The main theorem shows   that if $ f_{1j}=T_j(f_1) $,  then necessarily any  weak solution $(v_{1j})_j$ associated   to $f_{1j}=T_j(f_1)\in L^\infty (\O)$ remains in a bounded set of $W_0^{1,\qp}(\O), $ and there exists a subsequence $\big(v_{1\s(j)}\big)_j$ associated to $f_{1j}\in L^\infty(\O)$ and a function $v\in W_0^{1,\qp}(\O)$ which satisfy $\nabla v_{1\s(j)}\to \nabla v$ and $v_{1\s(j)}\to v$ a.e in $\O$ and $v_{1j}$ .\\
Let us show that we have necessarily $\Nu_1\equiv\nabla v$.\\
Indeed, for $k>0,\ B=\tan^{-1}(T_k)$ is in $W^{1,\infty}(\R),\ B'(0)=0$ if $|\s|>k$. Then, according to Lemma \ref{lem730124},  $\f=v_{1\s(j)}$ and $B(u-\f) $ are suitable test functions for both equations (weak formulation and entropic-renormalized formulation), hence we then have after letting $k\to+\infty $:
$$\int_\O\Big[\WH a_\pp(\Nu_1)-\WH a_\pp(v_{1\s(j)})\Big]\cdot\dfrac{\nabla (u-v_{1\s(j)})dx}{1+|u_1-v_{1\s(j)}|^2}\LEQ\dfrac\pi2 \int_\O |f_1-T_{\s(j)}f_1|.$$
Letting $j\to+\infty$
$$ \int_\O \dfrac{\Delta (u_1;v)}{1+|u_1-v|^2}dx=0$$
from which $\Delta (u_1;v)=0$ a.e., so that $\nabla u_1=\nabla v.$\\
This result shows that the whole sequence $(v_j)$ must satisfy   $\DST\lim_{j\to+\infty}\int_\O|\Nu_1-\nabla v_j|dx=0.$\\
This remark shows us if $f_1$ and $f_2$ are in $L^1(\O)$, then  we have a subsequence $T_{\s(j)}(f_1)$,\ \ $T_{\s(j)}(f_2)$ whose weak solutions $\Big(v_{1\s(j)}\Big)_j$,\ \ $(v_{2\s(j)})_j$ satisfy
$$\lim_{j\to+\infty}\nabla v_{i\s(j)}(x)=\nabla v_i(x)\ \aein \O.$$
As before, we easily have
$$\int_\O\dfrac{\Delta[v_{1\s(j)};v_{2\s(j)}]}{1+|v_{1\s(j)}-v_{2\s(j)}|^2}dx\LEQ\dfrac\pi2\int_\O\Big|T_{\s(j)}f_1-T_{s(j)}f_2\Big|.$$
Letting $j\to+\infty$, we get
$$\int_\O\dfrac{\Delta[u_1;u_2]}{1+|u_1-u_2|^2}\LEQ\dfrac\pi 2\int|f_1-f_2|dx,$$
from which we get the uniqueness.
\HF
{\bf Acknowledgments}\\
  A. Gogatishvili started work on this project during his visit to J.M. Rakotoson to the Laboratory of Mathematics of the University of Poitiers, France, in April 2022. He thanks for his generous hospithetality and helpful atmosphere during the visit.  He has been partially supported by   the   Czech Academy of Sciences (RVO 67985840), by  the Czech Science Foundation (GAČR), grant no: 23-04720S  by the Shota Rustaveli National Science Foundation (SRNSF), grant no:  FR-21-12353, and by the grant Ministry of
  Science and Higher Education of the Republic of Kazakhstan (project no. AP14869887).
  \\

M. R. Formica is partially supported by University of Naples "Parthenope", Dept. of Economic and Legal Studies, project CoRNDiS, DM MUR 737/2021, CUP I55F21003620001. \\
M. R. Formica is member of Gruppo Nazionale per l'Analisi Matematica, la Probabilit\'a e le loro Applicazioni (GNAMPA) of the Istituto Nazionale di Alta Matematica (INdAM) and member of the UMI group "Teoria dell'Approssimazione e Applicazioni (T.A.A.)".
\\

The authors thank the anonymous referee for the careful reading of the paper.


\end{document}